\documentclass{article}
\usepackage{amsmath}
\usepackage{amsfonts}
\usepackage{xcolor}

\usepackage{graphicx} 
\usepackage{subcaption} 

\newcommand\cO{\mathcal{O}}

\newtheorem{theorem}{Theorem}

\begin{document}

\centerline{\bf \Large Couette-Taylor instabilities in the small gap regime:}

\bigskip

\centerline{\bf \Large the very counter-rotating case}

\bigskip

\centerline{Dongfen Bian\footnote{School of Mathematics and Statistics, Beijing Institute of Technology, Beijing 100081, P.R. China. Email: biandongfen@bit.edu.cn}, 
%E. Grenier\footnote{Academy of Mathematics and Systems Science, Chinese Academy of Sciences, 100190, Beijing, P.R.China, emmanuelgrenier@amss.ac.cn},
Gérard Iooss\footnote{Laboratoire J.A.Dieudonné, I.U.F., Université Côte d’Azur, Parc Valrose, 06108 Nice Cedex 02, France}
%Zhuolun Yang\footnote{Department of Mathematics, The Ohio State University, Columbus, OH 43210, USA. Email: yang.8242@osu.edu}
}

\begin{abstract}
In this paper, we study the Couette-Taylor instability of a viscous fluid between two rotating cylinders in the small‑gap, slow rescaled rotation rate, high Reynolds number regime, focusing on the very counter‑rotating case $\mu < \mu_c \approx -0.8$ where the primary instability is non‑axisymmetric. Starting from the Navier-Stokes equations, we derive a limit system that captures the leading‑order dynamics and compute the critical Taylor number $T_c(\mu)$ together with the critical axial and azimuthal wavenumbers. Near criticality, the weakly nonlinear behaviour is governed by a system of two coupled complex Ginzburg-Landau equations. All coefficients of this amplitude system including the cubic nonlinear terms are evaluated numerically from the linearised eigenfunctions and the associated adjoint problem. The reduced equations admit helico\"{i}dal waves (travelling in both the axial and azimuthal directions) and ribbon waves (standing axially, travelling azimuthally), and their existence and stability criteria are discussed. We also examine more exotic spatially modulated solutions that satisfy a third‑order dynamical system, whose complete classification remains an open challenge. 
\end{abstract}

%%%%%%%%%%%%%%%%%%%%%%%%%%%%%%%%%%%%%%%%%%%%%%%%%%%%%%%%%%%%%%%%

\section{Introduction}

%%%%%%%%%%%%%%%%%%%%%%%%%%%%%%%%%%%%%%%%%%%%%%%%%%%%%%%%%%%%%%%%

In this paper, we study the classical problem of the flow of a viscous fluid between two rotating cylinders,
motivated in particular by the works of M. Nagata \cite{Nagata86}, \cite{Nagata23}, \cite{Nagata24}
and we consider the Navier-Stokes equations
\begin{equation} \label{NS1}
\partial_t u + (u \cdot \nabla) u - \nu \Delta u + \nabla p = 0,
\end{equation}
\begin{equation} \label{NS2}
\nabla \cdot u = 0
\end{equation}
in the domain 
$$
\Omega = \Bigl\{ (r,\varphi,z) \quad  | \quad 
r_i < r < r_o \Bigr\}.
$$
In these equations, $u(t,x,y,z)$ denotes the velocity of the fluid, $p(t,x,y,z)$ its pressure,
and $\nu > 0$ the viscosity. Moreover, $(r,\varphi,z)$ are the cylindrical coordinates.
The domain $\Omega$ is the interval between two cylinders with the same axis at $r = 0$,
the inner one being of radius $r_i$ and the outer one being of radius $r_o$.

The inner cylinder (respectively the outer) rotates with an angular velocity
$\omega_i$ (respectively, $\omega_o$). We assume that the fluid ``sticks" to the boundary, namely that the velocity of the fluid equals the velocity
of the cylinders at $r = r_i$ and $r = r_o$.

The Couette-Taylor flow between independently rotating cylinders has long served as a canonical system for studying hydrodynamic stability and pattern formation, dating back to the classic experiments of Taylor~\cite{Taylor23} and the weakly nonlinear theories of Davey~\cite{Davey62} and DiPrima and Eagles~\cite{DiPrima77}.  In the narrow‑gap limit, asymptotic descriptions of the bifurcating Taylor vortices were developed for co‑rotating and weakly counter‑rotating cylinders by Nagata~\cite{Nagata86,Nagata23,Nagata24} and, in a spatial‑dynamics framework, in \cite{BGIY1,BGIY2}.  When the cylinders counter‑rotate strongly ($\mu < \mu_c \approx -0.8$), the primary instability becomes non‑axisymmetric, leading to travelling azimuthal waves, ribbons and more complex spatio‑temporal patterns~\cite{Iooss}.  Recent large‑scale numerical studies by Deguchi and Nagata~\cite{DeguchiNagata11} and Meseguer and Marques~\cite{Meseguer02} have charted the stability boundaries and secondary bifurcations in this regime, but a systematic analytic treatment that resolves both the axial and azimuthal modulations in the small‑gap limit has been missing.  The present paper fills this gap by deriving and analysing a system of two coupled Ginzburg-Landau equations that capture the leading‑order nonlinear dynamics near the non‑axisymmetric onset.

We define the gap $d$, the radii ratio $\eta$ and the ratio of the rotation rates $\mu$ by
$$
d = r_o - r_i, \quad \eta = \frac{r_i}{r_o}, \quad \mu = \frac{\omega_o}{\omega_i}.
$$
We work in a rotating frame, with a rotation rate $\Omega_{rf}$ defined by
$$
\Omega_{rf} = \frac{\omega_i + \omega_o}{2} =  \omega_i \frac{1+\mu}{2}
$$
and consider the stability of the Couette flow, defined by
\begin{equation} \label{Couette}
U(r) = r \, \Omega_{rot} (r) = (A - \Omega_{rf})r + \frac{B}{r},
\end{equation}
in this rotating frame, where $A$ and $B$ are constants such that the velocity of the fluid coincides with 
the velocity of the cylinders at the inner and outer cylinders, which leads to
$$
A = \omega_i \frac{\mu - \eta^2}{1 - \eta^2},
\quad 
B = \omega_i r_i^2 \frac{1 - \mu}{1 - \eta^2}.
$$
We focus on the small gap / slow rescaled rotation rate / high Reynolds regime (SGSRRHR) where
\begin{equation}\label{limit}
    \eta \to 1, 
\qquad
\widehat{\omega}=\omega_i \frac{d^2}{\nu}\to 0, 
\qquad
{\mathfrak R} = \frac{ \widehat \omega (1 - \mu)}{1 - \eta}
\to \infty ,
\end{equation}
and define the Taylor number to be
$$
T = 2 \widehat{\omega} \mathfrak{R}.
$$
The question of the linear stability of the Couette flow is very classical and it is well established \cite{Drazin}
that there exists a critical Taylor
number $T_c$ such that the Couette flow is linearly stable if $T < T_c$ and linearly unstable if $T > T_c$. The instability is called
the Couette-Taylor instability. In the particular case $\mu = 1$,
$T_c \approx 1708$. For $T > T_c$, close to $T_c$, the unstable modes are axisymmetric and correspond to Fourier wave numbers in the $z$
variable close to $\alpha_c \approx 3.1$.

In \cite{BGIY1} and \cite{BGIY2} we have studied in detail the dependency of $T_c$ on $\mu$ in the SGSRRHR limit with $T$ of order 1, and numerical evidence shows that

\begin{itemize}

\item If $\mu > \mu_c \approx -0.8$, $T_c$ is a decreasing function of $\mu$ and the first instability is axisymmetric.

\item If $\mu < \mu_c \approx -0.8$, the first instability is non-axisymmetric.

\end{itemize}

In \cite{BGIY1} and \cite{BGIY2}, respectively in the cases $\mu = 1$ and $-0.8 < \mu < 1$, we have described these instabilities
in the SGSRRHR regime and discussed how small amplitude and slowly varying 
solutions of the initial Navier-Stokes
equations can be approximated by solutions of a reduced Ginzburg-Landau equation. We have further described the time and space
periodic solutions of this Ginzburg-Landau equation. Using the theory of spatial dynamics, this leads to the existence of non-trivial
time and space periodic solutions of the full Navier-Stokes equations, which are, at leading order, described by this Ginzburg-Landau equation.

\medskip

In the present article, we focus on the case $\mu < \mu_c \approx -0.8$ which turns out to be more complex since the instabilities
are no longer axisymmetric, but are periodic in the orthoradial variable. 
We will show that the behavior of small amplitude and slowly varying solutions of the Navier-Stokes
equations is described not by one, but by two coupled Ginzburg-Landau equations and numerically compute the coefficients of this
reduced system.

% We will further describe non-trivial particular solutions of these two Ginzburg-Landau equations. These solutions which are helicoïdal or ribbon waves approximate true solutions of the Navier-Stokes
% equations. We also discuss the existence of more complex solutions, called exotic ones.

To this end, 
we follow the general strategy of spatial dynamics, as developed in particular in \cite{Iooss2}.
The aim is to construct time and space periodic solutions to the full Navier-Stokes system (\ref{NS1})-(\ref{NS2}) 
as small perturbations of time and space periodic solutions of a reduced nonlinear system which catches the main dynamics near a bifurcation
point. In our case, this system consists of two coupled Ginzburg-Landau equations.

The strategy may be split into the following three steps:

\begin{itemize}

\item First, derivation of a reduced system which describes the long wave and small amplitude perturbations of the constant state
near the bifurcation point.

\item Second, construction of non-trivial periodic solutions for this reduced system.

\item Third, construction of genuine solutions of the Navier-Stokes system starting from solutions of the reduced system.

\end{itemize}

This last step can be done by repeating the proofs in \cite{Iooss2}, up to minor changes.
We will omit this step here and refer to \cite{Iooss2} for more details.

In the current paper, we focus on the first two points and describe new classes of solutions of the reduced system which 
lead to new classes of solutions of the genuine Navier-Stokes equations.

Note that the case $\mu = -1$ is considered in \cite{Nagata24} but without the study of non-axisymmetric bifurcating solutions. 
Such non-axisymmetric bifurcating solutions were studied in \cite{Iooss} for the full Navier-Stokes equations, but not in the small gap limit, 
which is our aim here.

\medskip

We present three main results.  First, we derive the small‑gap limit system and perform 
a complete linear stability analysis, obtaining the critical Taylor number $T_c(\mu)$ 
and the corresponding non‑axisymmetric wavenumbers $\alpha_c(\mu)$, $\mathfrak{B}_c(\mu)$ 
for the strongly counter‑rotating regime $\mu<\mu_c$.  

Second, using a center‑manifold 
reduction, we derive a system of two coupled complex Ginzburg–Landau equations that 
govern the slow modulation of the left‑ and right‑travelling azimuthal waves near onset.  
More precisely, small amplitude solutions $U$ of the limit system are of the form
$$
U  =A(y,t)\zeta _{1}+B(y,t)\zeta _{2}+\overline{A}(y,t)\overline{\zeta _{1}}%
+\overline{B}(y,t)\overline{\zeta _{2}}.
$$
where, up to higher order terms,
\begin{eqnarray}
\partial _{t}A &=&(i\omega _{0}+a_{3}\tau )A+b_{1}\mathfrak{R}\partial
_{y}A-a_{4}\mathfrak{R}^{2}\partial _{y}^{2}A+A(b|A|^{2}+c|B|^{2}),
\label{G1} \\
\partial _{t}B &=&(i\omega _{0}+a_{3}\tau )B+b_{1}\mathfrak{R}\partial
_{y}B-a_{4}\mathfrak{R}^{2}\partial _{y}^{2}B+B(c|A|^{2}+b|B|^{2}). \label{G2}
\end{eqnarray}%
We compute all coefficients of this amplitude system, namely, the linear dispersion terms 
$a_1,a_2,a_3,a_4,a_5,b_1$ and the cubic nonlinearities $b$ and $c$, directly from the 
linearised eigenfunctions and the adjoint problem by spectral Chebyshev collocation.

Third, we show that the Ginzburg-Landau systems admit particular solutions
which are helicoidal and ribbon waves. More precisely, we have the following theorem.

\begin{theorem}
The system (\ref{G1})-(\ref{G2}) admits solutions of the form
\begin{align*}
A(t,y) &= r_0 e^{(i \omega_0 + i a_{3i} \tau + i a_{4i} \beta^2) t + i \beta (y + b_{1r} t) + i \theta_0(t)},
\\
B(t,y) &= r_1 e^{(i \omega_0 + i a_{3i} \tau + i a_{4i} \beta^2) t + i \beta (y + b_{1r} t) + i \theta_1(t)}
\end{align*}
{\color{black}
where $r_0$ and $r_1$ arer real numbers,
 $\theta_0(t)$ and $\theta_1(t)$ are linear in $t$ and where $\tau = T - T_c$. The ``i" and ``r" indices
denote the imaginary and real parts of a complex number.
When $r_0 \ne 0$ and $r_1 = 0$, the solutions are helicoïdal waves, and when $r_0 = r_1 \ne 0$, they are ribbon waves.
}
\end{theorem}

We refer to the proof of this Theorem for the explicit expressions of $r_0$, $r_1$, $\theta_0$
and $\theta_1$.
We will also discuss more complex solutions and discuss their stability.
The complete classification of the solutions of (\ref{G1})-(\ref{G2}) remains an open problem.

\medskip

The following sections are organized as follows. In Section 2,  we present the small‑gap limit and linear analysis. In Section 3, we derive the Ginzburg–Landau system and its particular solutions, and in the Appendices, we detail the computations of the coefficients.

%%%%%%%%%%%%%%%%%%%%%%%%%%%%%%%%%%%%%%%%%%%%%%%%%%%%%%%%%%%%%%%%%

\section{The small gap limit}

%%%%%%%%%%%%%%%%%%%%%%%%%%%%%%%%%%%%%%%%%%%%%%%%%%%%%%%%%%%%%%%%%

%%%%%%%%%%%%%%%%%%%%%%%%%%%%%%%%%%%%%%%%%%%%%%%%

\subsection{The limit system}

%%%%%%%%%%%%%%%%%%%%%%%%%%%%%%%%%%%%%%%%%%%%%%%%

As shown in \cite{BGIY1}, \cite{BGIY2} and \cite{Nagata23}, \cite{Nagata24}, in the SGSRRHR limit, 
the initial Navier-Stokes equations (\ref{NS1})-(\ref{NS2}) reduce to the following system that we call the ``limit Navier-Stokes" system
\begin{equation}
\partial _{t}U=\Delta _{\bot }U+x\mathfrak{R}\partial _{y}U-\nabla _{\bot
}p+\left( 
\begin{array}{c}
Tg(x)\hat{u}_{y} \\ 
0 \\ 
u_{x}%
\end{array}%
\right) -(U_{\bot }\cdot \nabla _{\bot })U+\frac{T(1-\mu )}{2}\left( 
\begin{array}{c}
\hat{u}_{y}^{2} \\ 
0 \\ 
0%
\end{array}%
\right) ,  \label{NSlimit}
\end{equation}%
\begin{equation*}
\nabla _{\bot }\cdot U_{\bot }+\mathfrak{R}\partial _{y}\hat{u}_{y}=0,
\end{equation*}%
posed in the domain $-1/2 \le x \le 1/2$, $y \in \mathbb{R}$ and $z \in \mathbb{R}$.
In these equations, $U$ is the velocity of the fluid in the rotating frame, the index $\bot$ denotes the $(x,z)$ components,%
\begin{equation*}
U=(U_{\bot },\hat{u}_{y}), \qquad U_{\bot }=(u_{x},u_{z}),
\qquad 
u_{y}=\mathfrak{R}\hat{u}_{y},
\end{equation*}%
the parameters $T$ and $\mathfrak{R}$ are defined by
\begin{eqnarray*}
T &=&\frac{2\omega _{i}^{2}r_{0}^{4}(1-\eta )^{3}(1-\mu )}{\nu ^{2}}, \\
\mathfrak{R} &=&\frac{\omega _{i}r_{0}^{2}(1-\eta )(1-\mu )}{\nu }
\end{eqnarray*}%'
and
\begin{equation*}
g(x)=\frac{1+\mu }{2}-(1-\mu )x,
\end{equation*}%
where $x$ is the rescaled distance to the cylinders.

The boundary conditions are 
\begin{equation*}
U=0,x=\pm 1/2.
\end{equation*}%
We note that $T$ and $\mathfrak{R}$ are such that, $\mu \ne 1$ being fixed,
\begin{equation*}
(1-\eta )^{3} = \cO \Bigl( \frac{\nu }{\omega _{i}r_{0}^{2}} \Bigr)^{2}.
\end{equation*}%
In particular, $T = \cO(1)$ as $\eta \rightarrow 1,$ while $\mathfrak{R} = \cO[(1-\eta)^{-1/2}]$ as in \cite{Nagata24}. 

The derivation of (\ref{NSlimit}) starting from (\ref{NS1})-(\ref{NS2}) is detailed in \cite{BGIY1} and \cite{BGIY2}
and its main steps are recalled in the Appendix  \ref{derivation} for the reader's convenience.

We note that we keep the parameter $\mathfrak{R}$ in this limit system, whereas $\mathfrak{R} \to \infty$ since
we will consider solutions $U$ which are slowly varying in $y$, with scales of variations of order $\mathfrak{R}^{-1}$ in this
variable. For such solutions, the operator $\mathfrak{R} \partial_y$ remains bounded as $\mathfrak{R}$ goes to infinity.
It is of course possible to remove the parameter from the limit system by rescaling $U$ in $y$, but the rescaled system
is less clear than (\ref{NSlimit}). Keeping $\mathfrak{R}$ also clearly shows the natural scales in the $y$ variable.

Note that $x$ is bounded and ranges between $-1/2$ and $+1/2$. On the contrary $y \in \mathbb{R}$.
Moreover, we focus on solutions which are periodic in $z$ with a fixed period $2\pi /\alpha _{c}$, where $\alpha_c$ will be defined below.

%%%%%%%%%%%%%%%%%%%%%%%%%%%%%%%%%%%%%

\subsection{Linearization}

%%%%%%%%%%%%%%%%%%%%%%%%%%%%%%%%%%%%%

As shown in \cite{BGIY2}, the linearization of (\ref{NSlimit}) for general perturbations of the form
\begin{equation*}
e^{i(\alpha z+\beta y)}\widehat{U}(x)
\end{equation*}%
leads to the eigenvalue problem (see \cite{BGIY2})
\begin{eqnarray*}
(\lambda +\alpha ^{2}-D^{2}-i\mathfrak{B}x)u_{x}+Dp &=&Tg(x)\widehat{u}_{y},
\\
(\lambda +\alpha ^{2}-D^{2}-i\mathfrak{B}x)\widehat{u}_{y} &=&u_{x}, \\
(\lambda +\alpha ^{2}-D^{2}-i\mathfrak{B}x)u_{z}+i\alpha p &=&0, \\
Du_{x}+i\mathfrak{B}\widehat{u}_{y}+i\alpha u_{z} &=&0,
\end{eqnarray*}%
where
\begin{equation*}
\mathfrak{B}=\beta \mathfrak{R}=O(1),
\end{equation*}%
and%
\begin{equation*}
u_{x}=\widehat{u}_{y}=u_{z}=0 \quad \hbox{when} \quad x=\pm 1/2.
\end{equation*}%
We recall that the function $g(x)$ depends on the parameter $\mu$, which is fixed.
This eigenvalue problem depends on the parameters $T$, $\mathfrak{B}$ and $\alpha$.
Numerically, we observe that there exists a critical Taylor number $T_1(\alpha,\mathfrak{B})$ such that
all the eigenvalues of the linearised problem have negative real parts if $T < T_1(\alpha,\mathfrak{B})$
and such that there exists one eigenvalue with a positive real part if $T > T_1(\alpha,\mathfrak{B})$.

We define the critical Taylor number $T_c(\mu)$ by
$$
T_c(\mu) = \min_\alpha \min_{\mathfrak{B}} T_1(\alpha,\mathfrak{B}).
$$
Numerically, we have observed in \cite{BGIY2}  that $T_c(\mu)$ is reached for $\mathfrak{B} = 0$ provided $\mu > \mu_c$
where $\mu_c \approx - 0.8$. In this regime, the first instabilities are axisymmetric. We refer to \cite{BGIY1} and
\cite{BGIY2} for a detailed analysis of this regime.

When $\mu < \mu_c \approx -0.8$, $T_c(\mu)$ is reached for some $\alpha_c(\mu) \ne 0$ and some $\mathfrak{B}_c(\mu) > 0$.
It is also reached by symmetry when $\mathfrak{B} = - \mathfrak{B}_c(\mu)$.
We denote by $\lambda_0$ the corresponding eigenvalue of the linearised variable, which depends on 
the bifurcation parameter
$$
\tau = T - T_c
$$
and on $\alpha$ and $\mathfrak{B}$. By definition of $T_c$, at criticality, $\lambda = i \omega_0$ where $\omega_0 \in \mathbb{R}$.
Moreover, $\bar \lambda_0 = - i \omega_0$ is also an eigenvalue.

We now expand $\lambda_0$ near $T = T_c$ (namely near $\tau = 0$), $\alpha = \alpha_c(\mu)$ and $\mathfrak{B} = \mathfrak{B}_c(\mu)$.
By definition, at criticality, $\lambda_0$ is of the form
\begin{eqnarray}
\lambda _{0} &=&i\omega _{0}+ia_{1}(\alpha ^{2}-\alpha
_{c}^{2})+a_{2}(\alpha ^{2}-\alpha _{c}^{2})^{2}+a_{3}\tau +ib_{1}(\mathfrak{%
B}-\mathfrak{B}_{c})+a_{4}(\mathfrak{B}-\mathfrak{B}_{c})^{2}  \notag \\
&&+2a_{5}(\mathfrak{B}-\mathfrak{B}_{c})(\alpha ^{2}-\alpha _{c}^{2})+O(\tau
(\alpha ^{2}-\alpha _{c}^{2})) + \cO(\tau (\mathfrak{B}-\mathfrak{B}_{c})),
\label{lambda0}
\end{eqnarray}%
where%
\begin{equation*}
\omega _{0},a_{1},b_{1}\in \mathbb{R},
\end{equation*}%
and%
\begin{equation*}
a_{3r}>0,
\end{equation*}%
since $\Re \lambda _{0}>0$ for $\tau >0.$ All the coefficients of (\ref{lambda0}) depend on $\mu$. 

Let us now study the manifold $\Re \lambda_0 = 0$. If $\alpha$ and $\mathfrak{B}$ are given, $\Re \lambda_0 = 0$ provided
\begin{eqnarray*}
&&
a_{3r}\tau +a_{2r}(\alpha ^{2}-\alpha _{c}^{2})^{2}+a_{4r}(\mathfrak{B}- \mathfrak{B}_{c})^{2}+2a_{5r}(\mathfrak{B}
-\mathfrak{B}_{c})(\alpha^{2}-\alpha _{c}^{2}) 
\\
&&
+ \cO(\tau (\alpha ^{2}-\alpha _{c}^{2})) + \cO(\tau (\mathfrak{B}-\mathfrak{B}%
_{c})) = 0.
\end{eqnarray*}%
By definition of $T_c$, this equation on $\tau$ has no negative solution. This leads to the constraints
$$
a_{2r} + a_{4r} < 0,
\qquad 
a_{2r}a_{4r} > a_{5r}^{2},
$$
which expresses that the quadratic form giving the principal part of $\tau $
near $(\alpha _{c},\mathfrak{B}_{c})$ is negative definite.

The computation of the coefficients $a_1$, $a_2$, $a_3$, $a_4$ and $a_5$ is detailed in Appendix $4.2$.
In the sequel we will consider only periodic solutions in $z$, with period $2 \pi / \alpha_c(\mu)$.
Then for $\tau > 0$ close to $0$ and for $\mathfrak{B}$ close to $\mathfrak{B}_c$, the eigenvalue with the largest
real part of the linearised system is $\lambda_0$ which depends on $\tau$ and $\mathfrak{B}$.
For $\mathfrak{B} \ne \mathfrak{B}_c$, $\Re \lambda_0 < 0$ if $\tau = 0$ and more generally if $\tau$ is small enough,
and becomes positive when $\tau$ becomes larger (depending on $\mathfrak{B}$).

We will denote by $\zeta_1$ the eigenvector corresponding to $\alpha_c > 0$ and by $\zeta_2$ the eigenvector corresponding
to $- \alpha_c$. They are of the form
\begin{align*}
\zeta _{1} & = e^{i(\alpha _{c}z+\beta _{c}y)}U_{1}(x),
\\
\zeta _{2} &= S\zeta_{1} = e^{i(-\alpha _{c}z+\beta _{c}y)}SU_{1}(x)
\end{align*}
where $S$ is the representation of the symmetry $z\mapsto -z$.

%%%%%%%%%%%%%%%%%%%%%%%%%%%%%%%%%%%%%%%%%%%%%%%%%%%%%%%%%%%%%%%%%

\section{The Ginzburg-Landau system}

%%%%%%%%%%%%%%%%%%%%%%%%%%%%%%%%%%%%%%%%%%%%%%%%%%%%%%%%%%%%%%%%%

In this section, we establish the Ginzburg-Landau system which describes small amplitude solutions of the
limit Navier-Stokes equation and then construct particular solutions of this system called helicoïdal waves or ribbon waves
as well as more exotic solutions.

%%%%%%%%%%%%%%%%%%%%%%%%%%%%%%%%%%%%%%%%%%%%%%%%%%%%%%%%%%

\subsection{General setting of Ginzburg-Landau system}

%%%%%%%%%%%%%%%%%%%%%%%%%%%%%%%%%%%%%%%%%%%%%%%%%%%%%%%%%%

We now study solutions $U$ of Navier-Stokes which have a small amplitude and are slowly varying in $y$.
At leading order, we look for $U$ under the form
$$
U  =A(y,t)\zeta _{1}+B(y,t)\zeta _{2}+\overline{A}(y,t)\overline{\zeta _{1}}%
+\overline{B}(y,t)\overline{\zeta _{2}}.
$$
The dynamics of the amplitudes $A(y,t)$ and $B(y,t)$ is described by a Ginzburg-Landau system 
of the form (see \cite{Iooss} Chap. VIII)%
\begin{eqnarray*}
\partial _{t}A &=&f(A,\overline{A},B,\overline{B},\partial _{y}), \\
\partial _{t}B &=&g(A,\overline{A},B,\overline{B},\partial _{y}),
\end{eqnarray*}%
for some functions $f$ and $g$ which commute with the following symmetries (for every $h,k\in \mathbb{R}
)$%
\begin{eqnarray*}
(A,\overline{A},B,\overline{B},\partial _{y}) &\mapsto &(Ae^{ih},\overline{A}%
e^{-ih},Be^{-ih},\overline{B}e^{ih},\partial _{y}), \\
(A,\overline{A},B,\overline{B},\partial _{y}) &\mapsto &(Ae^{ik},\overline{A}%
e^{-ik},Be^{ik},\overline{B}e^{-ik},\partial _{y}), \\
(A,\overline{A},B,\overline{B},\partial _{y}) &\mapsto &(B,\overline{B},A,%
\overline{A},\partial _{y}).
\end{eqnarray*}%
We also need to recover the eigenvalues $\lambda_{0}$ and $\overline{\lambda _{0}}$ in its linear part. 
This implies that the principal part of the Ginzburg-Landau system is%
\begin{eqnarray}
\partial _{t}A &=&(i\omega _{0}+a_{3}\tau )A+b_{1}\mathfrak{R}\partial
_{y}A-a_{4}\mathfrak{R}^{2}\partial _{y}^{2}A+A(b|A|^{2}+c|B|^{2}),
\label{GL syst} \\
\partial _{t}B &=&(i\omega _{0}+a_{3}\tau )B+b_{1}\mathfrak{R}\partial
_{y}B-a_{4}\mathfrak{R}^{2}\partial _{y}^{2}B+B(c|A|^{2}+b|B|^{2}),  \notag
\end{eqnarray}%
which generalizes the standard system of Ginzburg-Landau equations described in \cite{Iooss} Chap. III.
The change of variable%
\begin{eqnarray*}
y &=&\mathfrak{R}(\widetilde{y}-b_{1}t), \\
A(y,t) &=&\widetilde{A}(\widetilde{y},t),
\qquad 
B(y,t)=\widetilde{B}(\widetilde{y},t),
\end{eqnarray*}%
leads to the system%
\begin{eqnarray*}
\partial _{t}\widetilde{A} &=&(i\omega _{0}+a_{3}\tau )\widetilde{A}%
-a_{4}\partial _{\widetilde{y}}^{2}\widetilde{A}+\widetilde{A}(b|\widetilde{A%
}|^{2}+c|\widetilde{B}|^{2}), \\
\partial _{t}\widetilde{B} &=&(i\omega _{0}+a_{3}\tau )\widetilde{B}%
-a_{4}\partial _{\widetilde{y}}^{2}\widetilde{B}+\widetilde{B}(c|\widetilde{A%
}|^{2}+b|\widetilde{B}|^{2}),
\end{eqnarray*}%
with 
\begin{equation*}
a_{3r}>0,
\qquad 
a_{4r}<0.
\end{equation*}
The coefficients $b$ and $c$ are studied in Appendix $4.3$.

We first consider particular solutions of (\ref{GL syst}) which are periodic in $y$ and only depend on the $t$ variable,
of the form
\begin{equation*}
\widetilde{A}(\widetilde{y},t)=\widehat{A}(t)e^{i\beta \widetilde{y}},
\qquad
\widetilde{B}(\widetilde{y},t)=\widehat{B}(t)e^{i\beta \widetilde{y}}.
\end{equation*}%
For such solutions, (\ref{GL syst}) reduces to the following system of two complex ordinary differential equations
\begin{eqnarray*}
\partial _{t}\widehat{A} &=&(i\omega _{0}+a_{3}\tau +a_{4}\beta ^{2})%
\widehat{A}+\widehat{A}(b|\widehat{A}|^{2}+c|\widehat{B}|^{2}), \\
\partial _{t}\widehat{B} &=&(i\omega _{0}+a_{3}\tau +a_{4}\beta ^{2})%
\widehat{B}+\widehat{B}(c|\widehat{A}|^{2}+b|\widehat{B}|^{2}).
\end{eqnarray*}%
If moreover we introduce the modulus and argument of $\widehat A$ and $\widehat B$ and write these two functions under the form
\begin{eqnarray*}
\widehat{A}(t) &=&r_{0}(t)e^{i(\omega _{0}+a_{3i}\tau +a_{4i}\beta
^{2})t+i\theta _{0}(t)}, \\
\widehat{B}(t) &=&r_{1}(t)e^{i(\omega _{0}+a_{3i}\tau +a_{4i}\beta
^{2})t+i\theta _{1}(t)},
\end{eqnarray*}%
we obtain the following system of four real ordinary differential equations
\begin{eqnarray}
\partial _{t}r_{0} &=&(a_{3r}\tau +a_{4r}\beta
^{2})r_{0}+r_{0}(b_{r}r_{0}^{2}+c_{r}r_{1}^{2}), \label{r1} \\
\partial _{t}r_{1} &=&(a_{3r}\tau +a_{4r}\beta
^{2})r_{1}+r_{1}(c_{r}r_{0}^{2}+b_{r}r_{1}^{2}), \label{r2} \\
\partial _{t}\theta _{0} &=&b_{i}r_{0}^{2}+c_{i}r_{1}^{2}, \label{r3} \\
\partial _{t}\theta _{1} &=&c_{i}r_{0}^{2}+b_{i}r_{1}^{2}. \label{r4}
\end{eqnarray}
Note that the equations on $r_0$ and $r_1$ are decoupled from the equations on $\theta_0$ and $\theta_1$.
Moreover, the two last equations are easily solved once $r_0$ and $r_1$ are known.

We now study the fixed points of (\ref{r1})-(\ref{r2}), which leads to so called helicoïdal and ribbon waves.
Up to symmetries, the system (\ref{r1})-(\ref{r2}) has three fixed points:

\begin{itemize}
\item $r_0 = r_1 =0$ which leads to $U = 0$.

\item $r_1 = 0$ and
$$
r_0^2 = - {a_{3r} \tau + a_{4r} \beta^2 \over b_r}
$$
which leads to helicoïdal waves, described in the next section.

\item $r_0 = r_1$ which leads to ribbon waves, studied in Section \ref{ribbon}.

\end{itemize}

%%%%%%%%%%%%%%%%%%%%%%%%%%%%%%%%%%%%%%%%%

\subsection{Helico\"{\i}dal waves}

%%%%%%%%%%%%%%%%%%%%%%%%%%%%%%%%%%%%%%%%%

We start with the particular case where $r_1 = 0$ or $r_0 = 0$, following 
the discussion of \cite{Iooss} {\color{black} (in particular Chap. III, Section III.2.1, page $46$)}.
A first solution is obtained by choosing $r_1 = 0$ (namely $\widehat B = 0$), together with 
\begin{equation*}
r_{0}^{2}=-\frac{a_{3r}\tau +a_{4r}\beta ^{2}}{b_{r}}, 
\qquad 
\theta _{0} = -\frac{b_{i}(a_{3r}\tau +a_{4r}\beta ^{2})}{b_{r}}t+\theta _{0}^{\ast }.
\end{equation*}
There exists a similar symmetric solution, with this time $r_0 = 0$.
Such a solution is called an \emph{helico\"{\i}dal waves} since it travels in both the azimuthal and the axial directions.

It exists provided $a_{3r} \tau + a_{4r} \beta^2$ and $b_r$ have opposite signs, namely provided 
\begin{equation*}
(a_{3r}\tau +a_{4r}\beta ^{2})b_{r}<0.
\end{equation*}%
If $b_{r}<0$, as $a_{4r}< 0$, they exist only on the supercritical side and more precisely when
\begin{equation*}
\tau >\frac{(-a_{4r})}{a_{3r}}\beta ^{2}.
\end{equation*}%
If $b_{r}>0$ they exist when
\begin{equation*}
\tau <\frac{(-a_{4r})}{a_{3r}}\beta ^{2},
\end{equation*}%
namely on both sides of $\tau =0.$ 

The corresponding solutions for
Navier-Stokes equations take the form of travelling waves (for the symmetric
solution, change $\alpha _{c}$ in $-\alpha _{c}$ in the formula)%
\begin{equation*}
U(x,y,z,t)=V(x,\alpha _{c}z+\widetilde{\beta }y+\omega t),
\end{equation*}%
where $V$ is $2\pi $ - periodic in its second argument,  with a principal part
\begin{equation}
r_{0}e^{i(\alpha _{c}z+\widetilde{\beta }y+\omega t+\theta ^{\ast
})}U_{1}(x)+c.c.  \label{helicoidalwave}
\end{equation}%
with%
\begin{eqnarray*}
\omega  &=&\omega _{0}+b_{1}\beta +a_{3i}\tau +a_{4i}\beta
^{2}+b_{i}r_{0}^{2}+h.o.t., \\
\widetilde{\beta } &=&\beta _{c}+\beta /\mathfrak{R}.
\end{eqnarray*}%
It is possible to study the stability of these solutions with respect
to perturbations which are $2\pi /\alpha _{c}$ - periodic in $z,$ and $2\pi /\beta $ -
periodic in $y$ , as it is done for $\beta =0$ in \cite{Iooss}. The result is
that an helico\"{\i}dal wave (\ref{helicoidalwave}) is stable provided 
\begin{equation*}
b_{r}<0 \quad \hbox{and} \quad c_{r}-b_{r}<0,
\end{equation*}%
which means that these waves are stable only when they occur on the
supercritical side, with $\tau >|a_{4r}| \beta ^{2} / a_{3r}$.

\medskip

\begin{figure}[htbp]
    \centering
        \includegraphics[width=\textwidth]{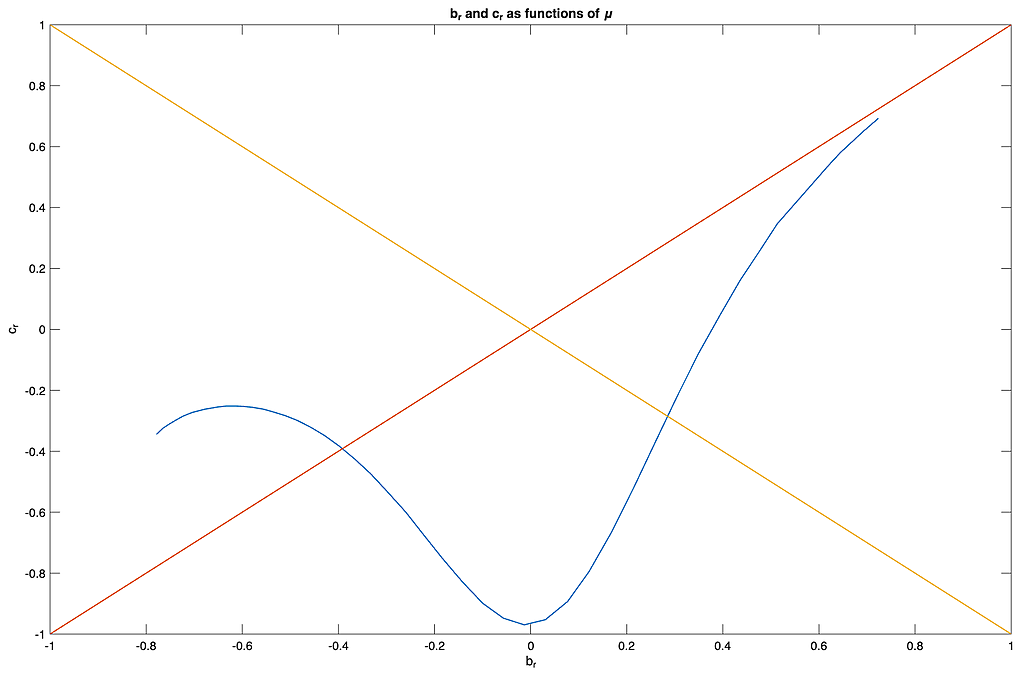}
        \caption{ {\color{black}
 Evolution of the point $(b_r,c_r)$ as a function of $\mu$ from $\mu = \mu_c$ to $\mu = -1$ 
        (from right to left) together with the diagonals $b_r = c_r$ and $b_r = - c_r$.
        The blue curve starts on the upper right part of the figure at $\mu = \mu_c$ and successively crosses
        $c_r = -b_r$ at $\mu = \mu_3$, $b_r =0$ at $\mu = \mu_1$ and $c_r = b_r$ at $\mu = \mu_2$. 
        It ends on the left part of the figure at $\mu = -1$.} }
        \label{figrebRec}
\end{figure}

{\color{black}
Figure \ref{figrebRec} shows $b_r$ and $c_r$ as a function of the parameter $\mu$
as $\mu$ goes from $\mu_c$ down to $-1$ (from right to left), together with the diagonals
$b_r = c_r$ and $b_r = - c_r$. We note that the curve $\mu \to (b_r(\mu),c_r(\mu))$ is not monotonic.

We recall that $\mu_c \approx -0.785$.
Numerically, we observe that

\begin{itemize}
    \item $b_r>0$ if $\mu > \mu_1 \approx -0.814$ and $b_r < 0$ if $\mu < \mu_1$.

\item $c_r - b_r > 0$ if $\mu < \mu_2 \approx - 0.848$ and $c_r - b_r < 0$ if $\mu > \mu_2$.

\item $c_r + b_r > 0$ if $\mu > \mu_3 \approx -0.8$ and $c_r + b_r < 0$ if $\mu < \mu_3$.

\end{itemize}
}

If $b_r > 0$, helico\"{i}dal waves exist on both sides, else they only exist on the supercritical side.
If $b_r < 0$ and $c_r -b_r <0$, they are stable.

Thus, if $\mu > \mu_1$, helicoïdal waves exist on both side. 
If $\mu_2 < \mu < \mu_1$, helicoïdal waves exist on the super-critical side and are stable.
If $\mu < \mu_2$, they exist on the super-critical side and are unstable.

%%%%%%%%%%%%%%%%%%%%%%%%%%%%%%%%%%%%%%%%%%%%%%%%%%%%%%%%%%%%%%%%

\subsection{Ribbon waves \label{ribbon}}

%%%%%%%%%%%%%%%%%%%%%%%%%%%%%%%%%%%%%%%%%%%%%%%%%%%%%%%%%%%%%%%%

The \emph{ribbon waves} correspond to 
\begin{equation*}
r_{0}^{2}=r_{1}^{2}=-\frac{a_{3r}\tau +a_{4r}\beta^{2}}{b_{r}+c_{r}},
\end{equation*}
$$
\theta _{k}=-\frac{(b_{i}+c_{i})(a_{3r}\tau +a_{4r}\beta^{2})}{%
b_{r}+c_{r}}t+\theta _{k}^{\ast }.
$$
They exist provided 
\begin{equation*}
(a_{3r}\tau +a_{4r}\beta ^{2})(b_{r}+c_{r})<0
\end{equation*}%
holds true. Hence, if $b_{r}+c_{r}<0$ they exist only on the supercritical side:%
\begin{equation*}
\tau >\frac{(-a_{4r})}{a_{3r}}\beta ^{2},
\end{equation*}%
while if $b_{r}+c_{r}>0$ they exist for%
\begin{equation*}
\tau <\frac{(-a_{4r})}{a_{3r}}\beta ^{2},
\end{equation*}%
corresponding to both sides of $\tau =0.$ 

The corresponding solutions of the Navier-Stokes equations take the form of standing waves in axial direction,
travelling in the azimuthal direction%
\begin{equation*}
U(x,y,z,t)=V(x,\alpha _{c}z,\widetilde{\beta }y+\omega t),
\end{equation*}%
where $V$ is $2\pi $ periodic in its second and third arguments,  with
a principal part 
\begin{equation}
r_{0}e^{i(\alpha _{c}z+\widetilde{\beta }y+\omega t+\theta _{0}^{\ast
})}U_{1}(x)+r_{0}e^{i(-\alpha _{c}z+\widetilde{\beta }y+\omega t+\theta
_{1}^{\ast })}SU_{1}(x)+c.c.  \label{ribbonwave}
\end{equation}%
with%
\begin{eqnarray*}
\omega &=&\omega _{0}+b_{1}\beta +a_{3i}\tau +a_{4i}\beta
^{2}+(b_{i}+c_{i})r_{0}^{2}+h.o.t. \\
\widetilde{\beta } &=&\beta _{c}+\beta /\mathfrak{R}.
\end{eqnarray*}%
It is possible to study the stability of these solutions with respect
to perturbations which are $2\pi /\alpha _{c}$ - periodic in $z,$ and $2\pi /\beta $ -
periodic in $y$ , as it is done for $\beta =0$ in \cite{Iooss}. The result is
that a ribbon wave (\ref{ribbonwave}) is stable provided that%
\begin{equation*}
b_{r}+c_{r}<0
\quad \hbox{and} \quad 
b_{r}-c_{r}<0,
\end{equation*}%
which means that this type of waves is stable only when they occur on the
supercritical side and $\tau >|a_{4r}|\beta ^{2} / a_{3r}$.

\medskip

Numerically, if $\mu > \mu_3$, they exist on both side. If $\mu_2 < \mu < \mu_3$, there
exist on the super-critical size and are unstable.
If $\mu < \mu_2$, they exist on the super-critical size and are stable.

%%%%%%%%%%%%%%%%%%%%%%%%%%%%%%%%%%%%%%%%%%%%%

\subsection{Exotic solutions}

%%%%%%%%%%%%%%%%%%%%%%%%%%%%%%%%%%%%%%%%%%%%%

We now look for more general solutions, of the form
\begin{eqnarray*}
\widetilde{A} &=&r_{0}(y,t)e^{i(\omega _{0}+a_{3i}\tau )t+\theta _{0}(y,t)},
\\
\widetilde{B} &=&r_{1}(y,t)e^{i(\omega _{0}+a_{3i}\tau )t+ \theta _{1}(y,t)},
\end{eqnarray*}%
which leads to a system of real partial differential equations of order 8:%
\begin{align*}
\partial _{t}r_{0} & = a_{3r}\tau r_{0}-a_{4r}[\partial _{\widetilde{y}%
}^{2}r_{0}-r_{0}(\partial _{\widetilde{y}}\theta _{0})^{2}]+a_{4i}(2\partial
_{\widetilde{y}}r_{0}\partial _{\widetilde{y}}\theta _{0}+r_{0}\partial _{%
\widetilde{y}}^{2}\theta _{0})+r_{0}(b_{r}r_{0}^{2}+c_{r}r_{1}^{2}), 
\\
\partial _{t}r_{1} & =a_{3r}\tau r_{1}-a_{4r}[\partial _{\widetilde{y}%
}^{2}r_{1}-r_{1}(\partial _{\widetilde{y}}\theta _{1})^{2}]+a_{4i}(2\partial
_{\widetilde{y}}r_{1}\partial _{\widetilde{y}}\theta _{1}+r_{1}\partial _{%
\widetilde{y}}^{2}\theta _{1})+r_{1}(b_{r}r_{1}^{2}+c_{r}r_{0}^{2}),
\\
\partial _{t}\theta _{0} &=-\frac{a_{4r}}{r_{0}^{2}}\partial _{\widetilde{y}%
}(r_{0}^{2}\partial _{\widetilde{y}}\theta _{0})-\frac{a_{4i}}{r_{0}}%
[\partial _{\widetilde{y}}^{2}r_{0}-r_{0}(\partial _{\widetilde{y}}\theta
_{0})^{2}]+b_{i}r_{0}^{2}+c_{i}r_{1}^{2}, 
\\
\partial _{t}\theta _{1} &=-\frac{a_{4r}}{r_{1}^{2}}\partial _{\widetilde{y}%
}(r_{1}^{2}\partial _{\widetilde{y}}\theta _{1})-\frac{a_{4i}}{r_{1}}%
[\partial _{\widetilde{y}}^{2}r_{1}-r_{1}(\partial _{\widetilde{y}}\theta
_{1})^{2}]+b_{i}r_{1}^{2}+c_{i}r_{0}^{2}.
\end{align*}
Note that this time these equations are no longer decoupled. To go further, 
we restrict our analysis to solutions such that $r_{0}$ or $r_{1}$ identically vanishes, 
and to solutions such that $r_{0}=r_{1},\theta _{0}=\theta _{1}+const,$ with in both cases $%
\partial _{t}r_{0}=\partial _{t}r_{1}=0$ and $\partial _{t}\theta _{0}=\partial
_{t}\theta _{1}=\kappa _{0},$ where $\kappa _{0}$ is constant.

Assuming $r_{1}=0,$ and defining a new function $K_{0}$ by
\begin{equation*}
r_{0}^{2}\theta _{0}^{\prime }=K_{0},
\end{equation*}%
the solutions satisfy the following third order ordinary differential equation
\begin{eqnarray*}
\frac{a_{4r}}{r_{0}} \Bigl[ r_{0}^{\prime \prime }-\frac{K_{0}^{2}}{r_{0}^{3}} \Bigr]
-\frac{a_{4i}}{r_{0}^{2}}K_{0}^{\prime } &=&a_{3r}\tau +b_{r}r_{0}^{2}, 
\\
\frac{a_{4i}}{r_{0}} \Bigl[ r_{0}^{\prime \prime }-\frac{K_{0}^{2}}{r_{0}^{3}} \Bigr]
+ \frac{a_{4r}}{r_{0}^{2}}K_{0}^{\prime } &=&b_{i}r_{0}^{2}-\kappa _{0}.
\end{eqnarray*}%
This leads to
\begin{eqnarray}
r_{0}^{\prime \prime }-\frac{K_{0}^{2}}{r_{0}^{3}} &=&\frac{r_{0}}{%
|a_{4}|^{2}} \Bigl[ a_{3r}a_{4r}\tau +(b_{r}a_{4r}+b_{i}a_{4i})r_{0}^{2}-\kappa_{0}a_{4i} \Bigr], 
\label{reducedsysthelicoidal} 
\\
K_{0}^{\prime } &=&\frac{r_{0}^{2}}{|a_{4}|^{2}} \Bigl[-a_{3r}a_{4i}\tau + (b_{i}a_{4r}-b_{r}a_{4i})r_{0}^{2}-\kappa _{0}a_{4r} \Bigr].  
\notag
\end{eqnarray}
If we look for solutions with $r_{0}=r_{1},$ $\theta _{0}=\theta
_{1}+const,$ $\partial _{t}\theta _{0}=\partial _{t}\theta _{1}=\kappa _{0}$
and defining $K_{0}$ as above, we obtain again the following third order ordinary differential equations 
\begin{eqnarray*}
\frac{a_{4r}}{r_{0}} \Bigl[r_{0}^{\prime \prime }-\frac{K_{0}^{2}}{r_{0}^{3}} \Bigr]-%
\frac{a_{4i}}{r_{0}^{2}}K_{0}^{\prime } &=&a_{3r}\tau
+(b_{r}+c_{r})r_{0}^{2}, \\
\frac{a_{4i}}{r_{0}} \Bigl[ r_{0}^{\prime \prime }-\frac{K_{0}^{2}}{r_{0}^{3}} \Bigr]+%
\frac{a_{4r}}{r_{0}^{2}}K_{0}^{\prime } &=&(b_{i}+c_{i})r_{0}^{2}-\kappa
_{0},
\end{eqnarray*}%
which leads to
\begin{eqnarray}
r_{0}^{\prime \prime }-\frac{K_{0}^{2}}{r_{0}^{3}} &=&\frac{r_{0}}{%
|a_{4}|^{2}} 
\Bigl\{a_{3r}a_{4r}\tau + \Bigl[ (b_{r}+c_{r})a_{4r}+(b_{i}+c_{i})a_{4i} \Bigr]r_{0}^{2}-\kappa _{0}a_{4i} \Bigr\},
\label{reducedsystribbon} \\
K_{0}^{\prime } &=&\frac{r_{0}^{2}}{|a_{4}|^{2}}
\Bigl\{-a_{3r}a_{4i}\tau + \Bigl[ (b_{i}+c_{i})a_{4r}-(b_{r}+c_{r})a_{4i} \Bigr]r_{0}^{2}-\kappa _{0}a_{4r}\Bigr\}. 
\notag
\end{eqnarray}
Both these third order systems take the form 
\begin{eqnarray}
r_{0}^{\prime \prime } &=&\frac{K_{0}^{2}}{r_{0}^{3}}+r_{0}(-a_{0}\tau
+2c_{0}r_{0}^{2}+d_{0}\kappa _{0}),  \label{3rd order system} \\
K_{0}^{\prime } &=&r_{0}^{2}(\alpha _{0}\tau +\gamma _{0}r_{0}^{2}+\delta
_{0}\kappa _{0}),  \notag
\end{eqnarray}%
where 
\begin{equation*}
a_{0}=-\frac{a_{3r}a_{4r}}{|a_{4}|^{2}}>0,
\quad 
\alpha _{0}=-\frac{a_{3r}a_{4i}}{|a_{4}|^{2}},
\quad 
d_{0}=\frac{-a_{4i}}{|a_{4}|^{2}},
\quad \delta _{0}=\frac{-a_{4r}}{%
|a_{4}|^{2}},
\end{equation*}%
and for (\ref{reducedsysthelicoidal}) 
\begin{equation*}
2c_{0}=\frac{b_{r}a_{4r}+b_{i}a_{4i}}{|a_{4}|^{2}},
\quad 
\gamma _{0}=\frac{b_{i}a_{4r}-b_{r}a_{4i}}{|a_{4}|^{2}},
\end{equation*}%
while for (\ref{reducedsystribbon})%
\begin{equation*}
2c_{0}=\frac{(b_{r}+c_{r})a_{4r}+(b_{i}+c_{i})a_{4i}}{|a_{4}|^{2}},
\quad 
\gamma_{0}=\frac{(b_{i}+c_{i})a_{4r}-(b_{r}+c_{r})a_{4i}}{|a_{4}|^{2}}.
\end{equation*}
Looking for solutions of (\ref{3rd order system}) such that $K_{0}$ is
independent of $y$ gives $r_{0}$ independent of $y$, we come back to the
known solutions seen at previous section. The existence of other kind of solutions is open.

%%%%%%%%%%%%%%%%%%%%%%%%%%%%%%%%%%%%%%%%%%%%%%%%%%%%%%%%%%%%

\section{Appendix}

%%%%%%%%%%%%%%%%%%%%%%%%%%%%%%%%%%%%%%%%%%%%%%%%%%%%%%%%%%%%

%%%%%%%%%%%%%%%%%%%%%%%%%%%%%%%%%%%%%%%%%%%%%%%%%%%%%%%%%%%%

\subsection{Derivation of the limit Navier-Stokes equations \label{derivation}}

%%%%%%%%%%%%%%%%%%%%%%%%%%%%%%%%%%%%%%%%%%%%%%%%%%%%%%%%%%%%

In this section, we detail the computations leading to \eqref{NSlimit}.
Let $\Delta_{cyl}$ be the Laplace operator in
cylindrical coordinates, namely
$$
\Delta_{cyl} = \frac{\partial^2}{\partial r^2} + \frac{1}{r} \frac{\partial}{\partial r} 
+ \frac{1}{r^2} \frac{\partial^2}{\partial \varphi^2} + \frac{\partial^2}{\partial z^2}.
$$
We perform the change of variables
$$
\hat{x} = \frac{r - \bar{r}}{d}, \quad \hat{y} = \frac{\bar{r}}{d} \varphi, 
\quad \hat{z} = \frac{z}{d}, \quad \hat{t} =  \frac{\nu}{d^2} t ,
$$
where $\bar{r}$ is the average radius
$$
\bar{r} = \frac{r_o + r_i}{2} = \frac{1+ \eta}{2} r_o.
$$
We first expand $\Delta_{cyl} u$, $r^{-2} u$, $u_\varphi^2/r$, $U' + \Omega_{rot}$, $\Omega_{rot}$ in $\eta$ and $\mu$.

\begin{itemize}

\item Expansion of $\Delta_{cyl} u$: we have
$$
\Delta_{cyl} u = \Big(\frac{\nu}{d^3} \frac{\partial^2}{\partial \hat{x}^2}  
+ \frac{1}{\hat{x} d + \bar{r}} \frac{\nu}{d^2} \frac{\partial}{\partial \hat{x}} 
+ \frac{\nu}{d^3} \frac{\bar{r}^2}{(\hat{x} d + \bar{r})^2} \frac{\partial^2}{\partial \hat{y}^2} 
+ \frac{\nu}{d^3}\frac{\partial^2}{\partial \hat{z}^2} \Big) \hat{u}.
$$
We have
$$
\frac{d}{\hat{x} d + \bar{r}} = \frac{1-\eta}{\hat{x} (1 - \eta) + \frac{1+\eta}{2}} = \cO(1 - \eta),
$$
\begin{align*}
    \frac{\bar{r}^2}{(\hat{x} d + \bar{r})^2} 
    =& \frac{(\frac{1+\eta}{2})^2}{(\hat{x} (1 - \eta) + \frac{1+\eta}{2})^2}
    = 1 - \frac{\hat{x}(1-\eta) (\hat{x}(1-\eta) + (1+\eta))}{(\hat{x} (1 - \eta) + \frac{1+\eta}{2})^2}\\
    =& 1 + \cO(1 - \eta),
\end{align*}
therefore,
$$
\Delta_{cyl} u = \frac{\nu}{d^3} \Big( \frac{\partial^2}{\partial \hat{x}^2}  + \frac{\partial^2}{\partial \hat{y}^2} 
+\frac{\partial^2}{\partial \hat{z}^2} \Big) \hat{u} 
+ \cO(1 - \eta) \frac{\nu}{d^3} \Big( \frac{\partial \hat{u}}{\partial \hat{x}} + \frac{\partial^2 \hat{u}}{\partial \hat{y}^2}\Big).
$$

\item Expansion of $r^{-2} u$:

\begin{align*}
    \frac{u}{r^2} = \frac{\nu}{d^3} \frac{d^2}{(\hat{x} d + \bar{r})^2} \hat{u} = \cO(1-\eta)^2 \frac{\nu}{d^3} \hat{u}.
\end{align*}

\item Expansion of $u_\varphi^2/r$:
{\color{black}
Note that
\begin{align*}
\frac{d}{r} %&= \frac{d}{r_0} \frac{r_0}{r} 
%= (1- \eta) \frac{d \hat{x} + \bar{r}}{r_0}= (1-\eta)^2 \hat{x} + \frac{(1-\eta)(1+\eta)}{2}.
= (1 - \eta) + (1 - \eta)^2 \Bigl( {1 \over 2} - x \Bigr) + \cO( ( 1 - \eta)^3)
\end{align*}
}
and
$$
\frac{u_\varphi^2}{r} = \frac{\nu^2}{d^3} \frac{d}{r} \hat{u}_{\hat{y}}^2 .
%= \frac{\nu^2}{d^3} \Big( (1-\eta)^2 \hat{x} + \frac{(1-\eta)(1+\eta)}{2} \Big) \hat{u}_{\hat{y}}^2.
$$

\item Study of $U' + \Omega_{rot}$: 

{\begin{align*}
    U' + \Omega_{rot} =& 2 (A - \Omega_{rf})
    = 2 \omega_i \Big( \frac{\mu - \eta^2}{1- \eta^2} - \frac{1 + \mu}{2} \Big)
    = - \omega_i \frac{(1 - \mu)(1 + \eta^2)}{1 - \eta^2}.
\end{align*}}

\item Expansion of $\Omega_{rot}$:

\begin{align*}
     \Omega_{rot} =&  (A - \Omega_{rf}) + \frac{B}{r^2}
    = - \omega_i \frac{(1 - \mu)(1 + \eta^2)}{1 - \eta^2} + \omega_i \frac{1-\mu}{1-\eta^2} \frac{r_i^2}{r^2}\\
    =& - \omega_i \frac{1-\mu}{1-\eta^2} \Big( \frac{1 + \eta^2}{2} - \frac{r_i^2}{r^2} \Big).
\end{align*}
Note that
\begin{align*}
     \frac{r_i^2}{r^2} =& \frac{r_i^2}{(\hat{x} d + \bar{r})^2}
     =  \frac{\eta^2}{(\hat{x} (1 - \eta) + \frac{1+\eta}{2})^2}
     = \Big( \frac{\eta}{(\hat{x} + \frac{1}{2}) (1 - \eta) + \eta} \Big)^2\\
     =& \Big( 1 + \frac{1-\eta}{\eta} (\hat{x} + \frac{1}{2})  \Big)^{-2}\\
     =& 1 - 2  \Big(\frac{1-\eta}{\eta} (\hat{x} + \frac{1}{2})  \Big) + 3  \Big( \frac{1-\eta}{\eta} (\hat{x} + \frac{1}{2})  \Big)^2 + \cdots.
\end{align*}
Hence,
\begin{align*}
    &\frac{1}{1-\eta^2} \Big( \frac{1 + \eta^2}{2} - \frac{r_i^2}{r^2} \Big) \\
    =& \frac{1}{1-\eta^2} \Big[ \frac{\eta^2 - 1}{2} + 2  \Big(\frac{1-\eta}{\eta} (\hat{x} 
    + \frac{1}{2})  \Big) - 3  \Big( \frac{1-\eta}{\eta} (\hat{x} + \frac{1}{2})  \Big)^2 + \cdots \Big]\\
    =& - \frac{1}{2} + \frac{1}{\eta(1+\eta)}(2 \hat{x} + 1) + \cO(1 - \eta).
\end{align*}
As
$$\frac{2}{\eta(1+\eta)} - 1 = - \frac{(\eta+2)(\eta - 1)}{\eta (1+\eta)} = \cO(1-\eta),$$
we have
$$
\frac{1}{1-\eta^2} \Big( \frac{1 + \eta^2}{2} - \frac{r_i^2}{r^2} \Big) = \hat{x} + \cO(1 - \eta)
$$
and hence
$$
\Omega_{rot} = - \omega_i (1 - \mu)  \Bigl[ \hat{x} + \cO(1 - \eta)\Bigr].
$$
\end{itemize}

We now adimensionalize the Navier-Stokes equations written in cylindrical coordinates.
In cylindrical coordinates, a (not necessarily axisymmetric) perturbation $(u_r, u_\varphi, u_z, p)$ of the Couette flow \eqref{Couette} satisfies the system (see \cite{Nagata23})
\begin{align*}
   \Big( \frac{\partial}{\partial t} - \nu (\Delta_{cyl} - &\frac{1}{r^2}) \Big) u_r 
   + \nu \frac{2}{r^2} \frac{\partial u_\varphi}{\partial \varphi} 
   + \frac{\partial p}{\partial r} 
   \\
   &= - \Omega_{rot} \frac{\partial u_r}{\partial \varphi} + 2\Omega_{rot} u_\varphi  
   + \frac{u_\varphi^2}{r} + 2\Omega_{rf} u_\varphi - (u \cdot \nabla) u_r,
   \\
   \Big( \frac{\partial}{\partial t} - \nu (\Delta_{cyl} &- \frac{1}{r^2}) \Big) u_\varphi 
   - \nu \frac{2}{r^2} \frac{\partial u_r}{\partial \varphi} + \frac{1}{r} \frac{\partial p}{\partial \varphi} 
   \\
   &= - \Omega_{rot} \frac{\partial u_\varphi}{\partial \varphi} { - U' u_r - \Omega_{rot} u_r} 
   - \frac{u_r u_\varphi}{r} - 2\Omega_{rf} u_r  - (u \cdot \nabla) u_\varphi,
   \\
\Big( \frac{\partial}{\partial t} - &\nu \Delta_{cyl}  \Big) u_z 
+ \frac{\partial p}{\partial z} 
= - \Omega_{rot} \frac{\partial u_z}{\partial \varphi}  - (u \cdot \nabla) u_z,
\end{align*}
together with the incompressibility condition, which in
cylindrical coordinates reads
$$
\frac{\partial u_r}{\partial r} + \frac{u_r}{r} + \frac{1}{r} \frac{\partial u_\varphi}{\partial \varphi} 
+ \frac{\partial u_z}{\partial z} = 0.
$$
Here the Coriolis term is added on the right-hand side.
We rescale the velocity and the pressure by
\begin{align*}
    u(t,r,\varphi,z) = \frac{\nu}{d} \hat{u}(\hat{t}, \hat{x}, \hat{y},\hat{z}) ,\\
    p(t,r,\varphi,z) =\frac{\nu^2}{d^2} \hat{p}(\hat{t}, \hat{x}, \hat{y},\hat{z}) .
\end{align*}
This leads to
\begin{align*}
   \Big( \frac{\partial}{\partial \hat{t}} - \hat\Delta \Big) \hat{u}_{\hat{x}} 
   &+ \frac{\partial \hat p}{\partial \hat{x}} 
   = \mathfrak{R}\hat{x} \frac{\partial \hat{u}_{\hat{x}}}{\partial \hat{y}} 
   + \hat{\omega}(1 + \mu) \hat{u}_{\hat{y}} - (\hat{u} \cdot \nabla) \hat{u}_{\hat{x}} 
   +\cO \Bigl(\hat{\omega} (1 - \mu)\frac{\partial \hat{u}_{\hat{x}}}{\partial \hat{y}} \Bigr)
  \\
   &+ \cO(1 - \eta)^2 \hat{u}_{\hat{x}} - 2\Big(\hat{\omega} (1 - \mu)  \hat{x} + \cO(1 - \eta) \Big)\hat{u}_{\hat{y}} 
   +  \cO(1 - \eta)  \frac{\partial \hat{u}_{\hat{y}}}{ \partial \hat{y}} \\
   &+ \Big((1 - \eta)^2 x + \frac{(1-\eta)(1+\eta)}{2} \Big) \hat{u}_{\hat{y}}^2,
 \\
   \Big( \frac{\partial}{\partial \hat{t}} - \hat\Delta \Big) \hat{u}_{\hat{y}} 
   &+ \frac{\partial \hat p}{\partial \hat{y}} 
= \mathfrak{R}\hat{x}\frac{\partial \hat{u}_{\hat{y}}}{\partial \hat{y}} + 
\hat{\omega} \frac{(1 - \mu)(1 + \eta^2)}{1 - \eta^2} \hat{u}_{\hat{x}} - \hat{\omega}(1 + \mu) \hat{u}_{\hat{x}} - (\hat{u} \cdot \nabla) \hat{u}_{\hat{y}}
  \\
   & 
+  \cO(1 - \eta) \Big( \hat{u}_{\hat{x}} \hat{u}_{\hat{y}} + \frac{\partial \hat{u}_{\hat{x}}}{ \partial \hat{y}} \Bigr) 
+ \cO \Bigl(\hat{\omega} (1 - \mu)\frac{\partial \hat{u}_{\hat{y}}}{\partial \hat{y}} \Bigr)+ \cO(1 - \eta)^2 \hat{u}_{\hat{y}},
\\
   \Big( \frac{\partial}{\partial \hat{t}} - \hat\Delta \Big) \hat{u}_{\hat{z}} 
   &+ \frac{\partial \hat p}{\partial \hat{z}} 
   = \mathfrak{R}\hat{x}\frac{\partial \hat{u}_{\hat{z}}}{\partial \hat{y}}  - (\hat{u} \cdot \nabla) \hat{u}_{\hat{z}} 
   +\cO \Bigl( \hat{\omega} (1 - \mu)\frac{\partial \hat{u}_{\hat{z}}}{\partial \hat{y}} \Bigr) ,
\end{align*}
$$
    \frac{\partial \hat{u}_{\hat{x}}}{\partial \hat{x}} 
    +\frac{\partial \hat{u}_{\hat{y}}}{\partial \hat{y}} + \frac{\partial \hat{u}_{\hat{z}}}{\partial \hat{z}} 
    = \cO(1-\eta) \Bigl( \hat{u}_{\hat{x}} + \frac{\partial \hat{u}_{\hat{y}}}{ \partial \hat{y}} \Bigr).\nonumber
$$
In this system, $\cO(\cdot)$ denotes various smooth functions that go to zero in the limit \eqref{limit}.
As $\eta$ goes to $1$, under the scaling (\ref{limit}), the limiting system is
\begin{eqnarray}
(\partial _{t}-\Delta _{\bot })u_{\bot }+\nabla _{\bot }p &=&\mathfrak{R}%
x\partial _{y}u_{\bot }+T g(x) (%
\widehat{u}_{y},0)^{t}\nonumber
\\&&-(u_{\bot }\cdot \nabla _{\bot })u_{\bot } 
+\frac{T}{2}(1-\mu)(\widehat{u}_y^2,0)^{t},\label{nl1}
\\
(\partial _{t}-\Delta _{\bot })\widehat{u}_{y} &=&\mathfrak{R}x\partial _{y}%
\widehat{u}_{y}+u_{x}-(u_{\bot }\cdot \nabla _{\bot })\widehat{u}_{y}, \label{nl2}
\\
\nabla _{\bot }\cdot u_{\bot }+\mathfrak{R}\partial _{y}\widehat{u}_{y} &=&0, \label{nl3}
\end{eqnarray}%
where the subscript $\bot $ denotes the  components in the $(x,z)$ plane, $\widehat{u}_{y}=\mathfrak{R}u_{y}$, and
$$
g(x)=\frac{1+\mu}{2}-(1-\mu)x.
$$
In this system, we keep the linear terms with $\mathfrak{R} \partial_y$ 
since this operator remains bounded in the limiting regime described previously.

%%%%%%%%%%%%%%%%%%%%%%%%%%%%%%%%%%%%%%%%%%%%%%%%%%%%%%%%%%%%%%%%%%%%%%%%%%%%%%%%%%%%%%%%%%%%%%

\subsection{Range of $(i\protect\omega _{0}-L_{0})$ and vector $\protect \zeta _{1}^{\ast }$ \label{appendixrange}}

%%%%%%%%%%%%%%%%%%%%%%%%%%%%%%%%%%%%%%%%%%%%%%%%%%%%%%%%%%%%%%%%%%%%%%%%%%%%%%%%%%%%%%%%%%%%%%

We define the scalar product $\langle U, V \rangle$ of two vector fields $U = (u_x,u_y,u_z)$ and $V = (v_x,v_y,v_z)$ by
\begin{equation}
\langle U, V \rangle = 
\int_{-1/2}^{1/2} \Bigl[ Du_x \overline{Dv_x} +\alpha _{c}^{2}u_x
\overline{v_x} + u_{y}\overline{v_y} \Bigr] \, dx.
\end{equation}%

We have to take care that the divergence free vector field $U$ is no longer orthogonal to
gradients. Let us solve the equation%
\begin{equation*}
(i\omega _{0}-L_{0})U=F
\end{equation*}%
where%
\begin{eqnarray*}
F &=&e^{i(\alpha _{c}z+\beta _{c}y)}\Phi (x),
\qquad 
\Phi =(\phi _{x},\phi _{y},\phi
_{z}), \\
0 &=&D\phi _{x}+i\mathfrak{B}\phi _{y}+i\alpha \phi _{z}.
\end{eqnarray*}%
The vector $\Phi $ satisfies the same boundary conditions as $U$ at $x=\pm 1/2.$ We
look for $\widehat{U}(x)$ such that%
\begin{eqnarray*}
U &=&e^{i(\alpha _{c}z+\beta _{c}y)}\widehat{U}(x),
\qquad
\widehat{U}(x)=(u_{x},u_{y},u_{z}), \\
0 &=&Du_{x}+i\mathfrak{B}_{c}u_{y}+i\alpha u_{z}.
\end{eqnarray*}%
This gives%
\begin{eqnarray*}
(i\omega _{0}+\alpha _{c}^{2}-D^{2}-i\mathfrak{B}%
_{c}x)u_{x}+Dp-T_{c}g(x)u_{y} &=&\phi _{x}, \\
(i\omega _{0}+\alpha _{c}^{2}-D^{2}-i\mathfrak{B}_{c}x)u_{y}-u_{x} &=&\phi
_{y}, \\
(i\omega _{0}+\alpha _{c}^{2}-D^{2}-i\mathfrak{B}_{c}x)u_{z}+i\alpha p
&=&\phi _{z}, \\
Du_{x}+i\mathfrak{B}_{c}u_{y}+i\alpha u_{z} &=&0,
\end{eqnarray*}%
which leads to%
\begin{eqnarray}
(i\omega _{0}+\alpha _{c}^{2}-D^{2}-i\mathfrak{B}_{c}x)(\alpha
_{c}^{2}-D^{2})u_{x}-\alpha _{c}^{2}T_{c}g(x)u_{y} &=&(\alpha
_{c}^{2}-D^{2})\phi _{x},  \label{linsystsingular} \\
(i\omega _{0}+\alpha _{c}^{2}-D^{2}-i\mathfrak{B}_{c}x)u_{y}-u_{x} &=&\phi
_{y}.  \notag
\end{eqnarray}%
We then define the eigenvector belonging to $i\omega _{0}$%
\begin{equation*}
\zeta _{1}=e^{i(\alpha _{c}z+\beta
_{c}y)}(u_{x}^{0},u_{y}^{0},u_{z}^{0})^{t},
\end{equation*}%
with%
\begin{eqnarray}
(i\omega _{0}+\alpha _{c}^{2}-D^{2}-i\mathfrak{B}_{c}x)(\alpha
_{c}^{2}-D^{2})u_{x}^{0}-\alpha _{c}^{2}T_{c}g(x)u_{y}^{0} &=&0,
\label{eigenvect} \\
(i\omega _{0}+\alpha _{c}^{2}-D^{2}-i\mathfrak{B_{c}}x)u_{y}^{0}-u_{x}^{0}
&=&0,  \notag
\end{eqnarray}%
and 
\begin{equation*}
u_{x}^{0}=Du_{x}^{0}=u_{y}^{0}=0,x=\pm 1/2.
\end{equation*}%
Let us look for \emph{the compatibility condition} for solving (\ref%
{linsystsingular}). We introduce the vector field%
\begin{eqnarray*}
\zeta _{1}^{\ast } &=&e^{i(\alpha _{c}z+\beta
_{c}y)}(v_{x}^{0},v_{y}^{0},v_{z}^{0}), \\
0 &=&Dv_{x}^{0}+i\mathfrak{B}_{c}v_{y}^{0}+i\alpha v_{z}^{0}, \\
v_{x}^{0} &=&Dv_{x}^{0}=v_{y}^{0}=0,x=\pm 1/2,
\end{eqnarray*}%
satisfying the adjoint system of (\ref{eigenvect})%
\begin{eqnarray}
(\alpha _{c}^{2}-D^{2})(-i\omega _{0}+\alpha _{c}^{2}-D^{2}+i\mathfrak{B}%
_{c}x)v_{x}^{0}-v_{y}^{0} &=&0,  \label{adj eigenvect} \\
(-i\omega _{0}+\alpha _{c}^{2}-D^{2}+i\mathfrak{B}_{c}x)v_{y}^{0}-\alpha
_{c}^{2}T_{c}g(x)v_{x}^{0} &=&0,  \notag
\end{eqnarray}%
then it is easy to check for (\ref{linsystsingular}) that $(\phi _{x},\phi
_{y})$ needs to satisfy%
\begin{equation*}
\int_{-1/2}^{1/2} \Bigl\{[(\alpha _{c}^{2}-D^{2})\phi _{x}]\overline{v_{x}^{0}}%
+\phi _{y}\overline{v_{y}^{0}} \Bigr\} \, dx=0,
\end{equation*}%
i.e.%
\begin{equation}
\int_{-1/2}^{1/2} \Bigl[ D\phi _{x}D\overline{v_{x}^{0}}+\alpha _{c}^{2}\phi _{x}%
\overline{v_{x}^{0}}+\phi _{y}\overline{v_{y}^{0}} \Bigr] \, dx=0,  \label{compatib}
\end{equation}%
which may be taken as a definition of the condition 
\begin{equation*}
\langle F,\zeta _{1}^{\ast }\rangle =0.
\end{equation*}%
We notice that (with the definition above)%
\begin{equation}
\langle \zeta _{1},\zeta _{1}^{\ast }\rangle 
=\int_{-1/2}^{1/2} \Bigl[ Du_{x}^{0}D%
\overline{v_{x}^{0}}+\alpha _{c}^{2}u_{x}^{0}\overline{v_{x}^{0}}+u_{y}^{0}%
\overline{v_{y}^{0}} \Bigr] \, dx  \label{scal prof dzeta1}
\end{equation}%
which should be non zero in the next calculations.

\begin{figure}[htbp]
    \centering
    \begin{subfigure}[b]{0.45\textwidth}
        \centering
        \includegraphics[width=\textwidth]{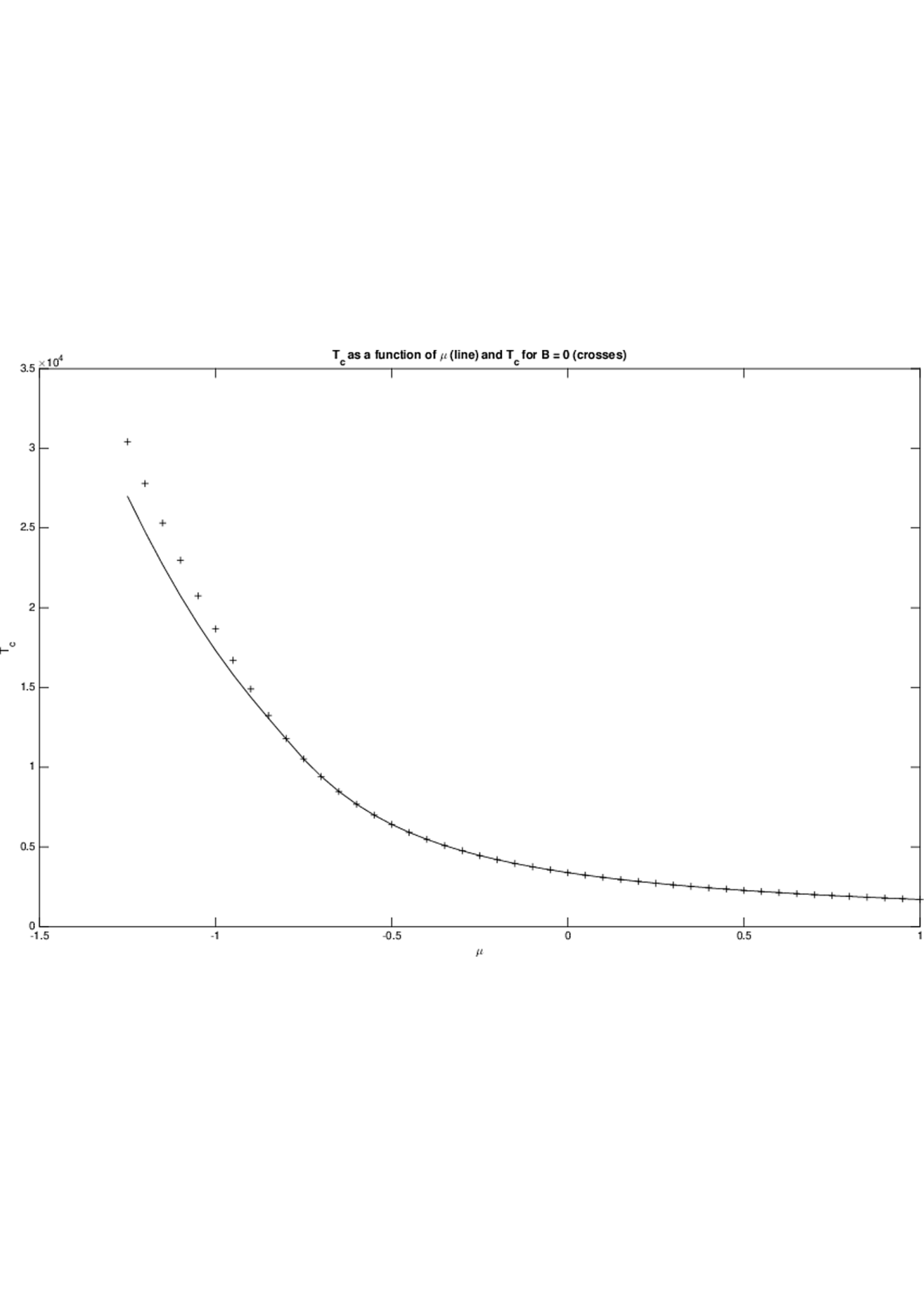}
        \caption{Critical Taylor number as a function of $\mu$ (solid line) and critical Taylor number
        corresponding to ${\mathfrak B} = 0$ }
        \label{fig:sub1}
    \end{subfigure}
    \hfill
    \begin{subfigure}[b]{0.45\textwidth}
        \centering
        \includegraphics[width=\textwidth]{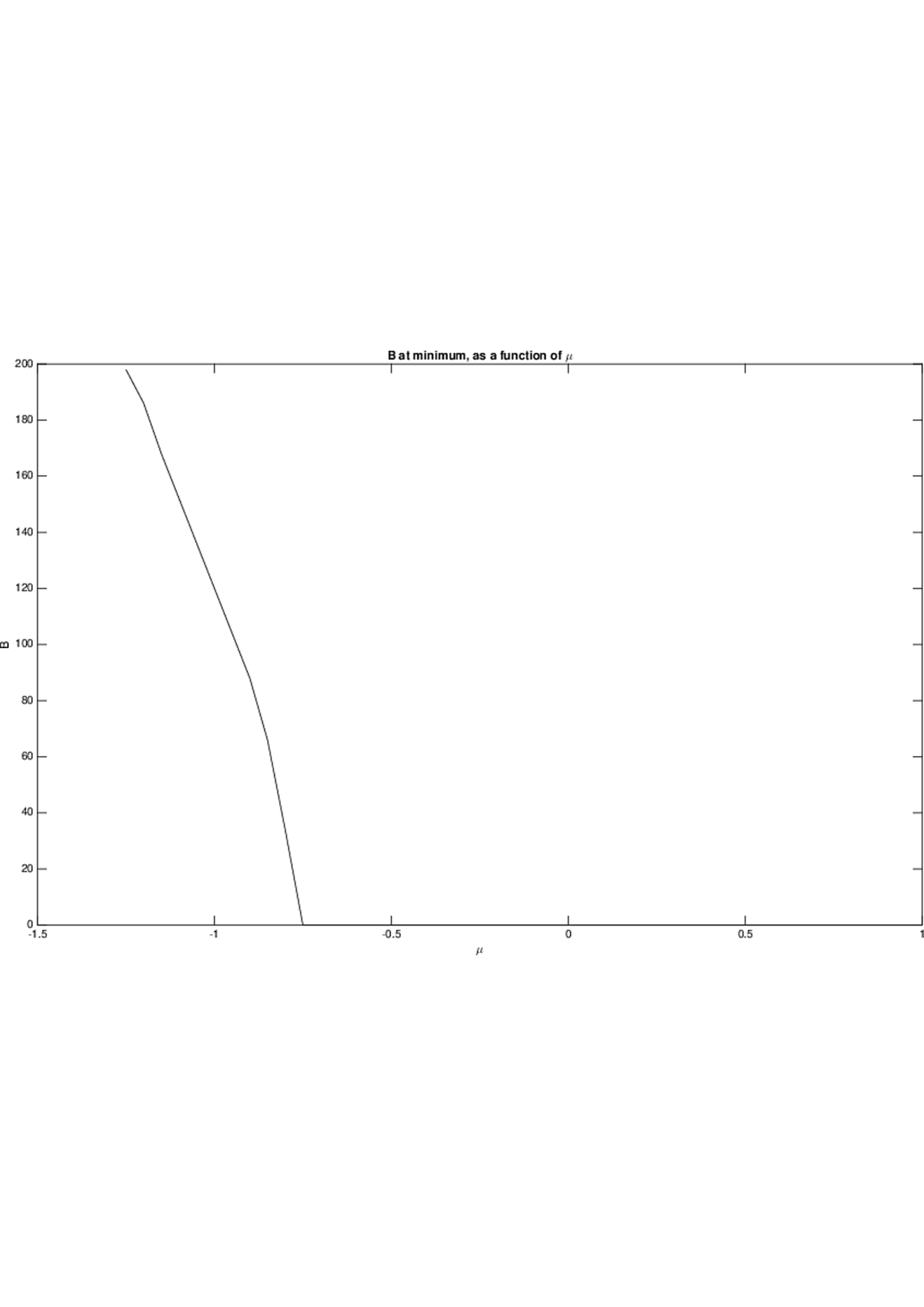}
        \caption{$\mathfrak{B}$ at the critical point as a function of $\mu$}
        \label{fig:sub2}
    \end{subfigure}
    \\[0.5cm] 
    \begin{subfigure}[b]{0.45\textwidth}
        \centering
        \includegraphics[width=\textwidth]{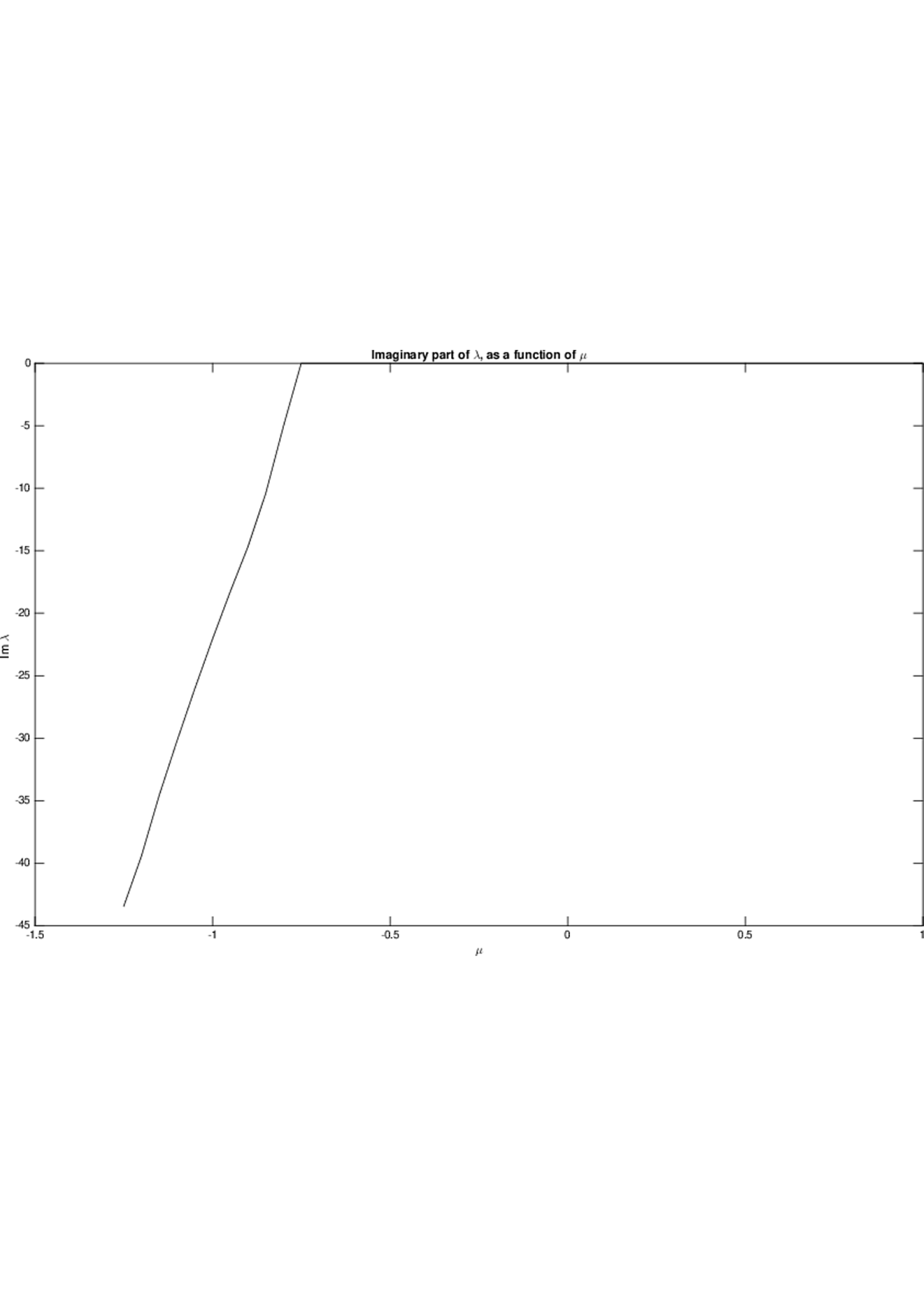}
        \caption{The imaginary part of the eigenvalue at the critical point as a function of $\mu$}
        \label{fig:sub3}
    \end{subfigure}
    \hfill
    \begin{subfigure}[b]{0.45\textwidth}
        \centering
        \includegraphics[width=\textwidth]{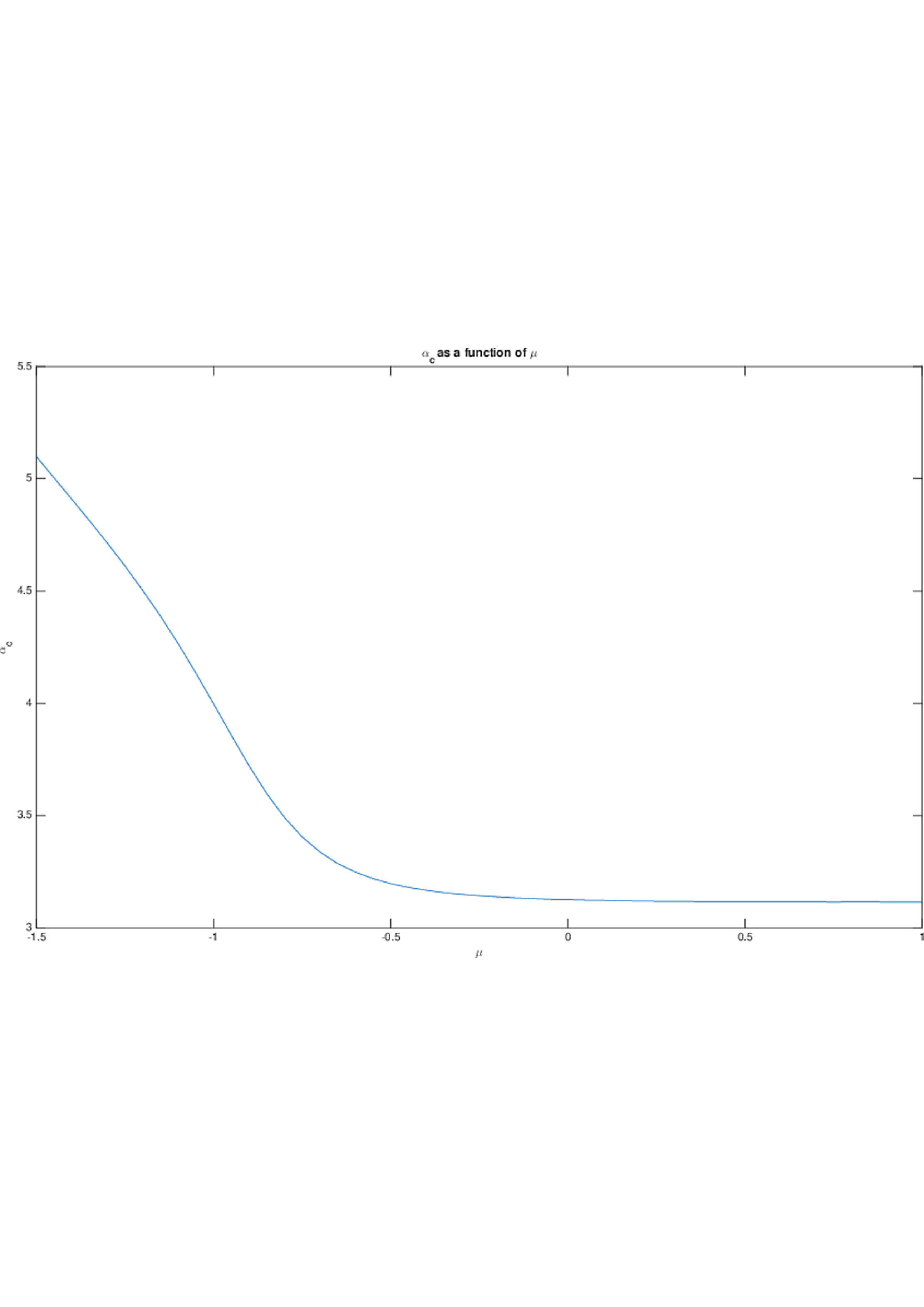}
        \caption{$\alpha$ at the critical point as a function of $\mu$}
        \label{fig:sub4}
    \end{subfigure}
    \caption{Numerical study of the critical point as a function of $\mu$. We observe that when $\mu < \mu_c \approx - 0,8$, the critical
    point is no longer axisymmetric ($\mathfrak{B} \ne 0$).}
    \label{fig:principale}
\end{figure}

%%%%%%%%%%%%%%%%%%%%%%%%%%%%%%%%%%%%%%%%%%%%%%%%%%%%%%%%%%%%%%%%%%%%%%%%%%%%%%%%%

\subsection{Computation of $a_{1},a_{3},b_{1},a_{4}$ \label{a134}}

%%%%%%%%%%%%%%%%%%%%%%%%%%%%%%%%%%%%%%%%%%%%%%%%%%%%%%%%%%%%%%%%%%%%%%%%%%%%%%%%%

To compute the coefficients $a_1$, $a_2$, $a_3$, $a_4$ and $a_5$, we rewrite
\begin{eqnarray}
(\lambda _{0}+\alpha ^{2}-D^{2}-i\mathfrak{B}x)(\alpha
^{2}-D^{2})u_{x}-\alpha ^{2}Tg(x)u_{y} &=&0,  \label{exp lambda0} \\
(\lambda _{0}+\alpha ^{2}-D^{2}-i\mathfrak{B}x)u_{y}-u_{x} &=&0,  \notag
\end{eqnarray}%
under the form
\begin{eqnarray}
\lambda _{0}(H_{0}+\widehat{\alpha }H_{010})U &=&L_{0}U
+ \widehat{\mathfrak{B}} L_{100}U + \widehat{\alpha }L_{010}U +\widehat T L_{001}U
\label{exp lambda b} \\
&&+\widehat{\alpha } \widehat{ \mathfrak{B} } L_{110}U + \widehat{\alpha } \widehat{T}  L_{011}U
+\widehat{\alpha }^{2}L_{020}U,  \notag
\end{eqnarray}%
where 
$$
\widehat T =T-T_{c},
\qquad
\widehat{\alpha }=\alpha ^{2}-\alpha _{c}^{2},
\qquad
\widehat{\mathfrak{B}} = \mathfrak{B} - \mathfrak{B_c}
$$
and where we have defined 
\begin{equation*}
H_{0}U=\left( 
\begin{array}{c}
(\alpha _{c}^{2}-D^{2})u_{x} \\ 
u_{y}%
\end{array}%
\right) ,
\qquad
H_{010}U=\left( 
\begin{array}{c}
u_{x} \\ 
0%
\end{array}%
\right) ,
\end{equation*}%
\begin{equation*}
L_{0}U=\left( 
\begin{array}{c}
- (\alpha _{c}^{2}-D^{2}-i\mathfrak{B}_{c}x)(\alpha_{c}^{2}-D^{2})u_{x} + \alpha_{c}^{2}T_{c}g(x)u_{y} 
\\ 
- (\alpha _{c}^{2}-D^{2}-i\mathfrak{B}_{c}x)u_{y} + u_{x}
\end{array}%
\right) ,
\end{equation*}%
\begin{equation*}
L_{001}U=\left( 
\begin{array}{c}
\alpha _{c}^{2}g(x)u_{y} \\ 
0%
\end{array}%
\right) ,
\quad
L_{100}U=\left( 
\begin{array}{c}
ix(\alpha _{c}^{2}-D^{2})u_{x} \\ 
ixu_{y}%
\end{array}%
\right) ,
\end{equation*}%
\begin{equation*}
L_{010}U=\left( 
\begin{array}{c}
-2(\alpha _{c}^{2}-D^{2})u_{x} + i {\mathfrak B}_c x u_x +T_{c}g(x)u_{y} \\ 
-u_{y}%
\end{array}%
\right) ,
\end{equation*}%
\begin{equation*}
L_{110}U=\left( 
\begin{array}{c}
ixu_{x} \\ 
0%
\end{array}%
\right) ,
\qquad
L_{011}U=\left( 
\begin{array}{c}
 g(x) u_{y} \\ 
0%
\end{array}%
\right) ,
\qquad 
L_{020}U=\left( 
\begin{array}{c}
-u_{x} \\ 
0%
\end{array}%
\right) .
\end{equation*}%
We look for an expansion of $U$ under the form
$$
U = U^0 + \widehat {\mathfrak B} U^{100} + \widehat \alpha U^{010} 
+ \widehat T U^{001}
+ \widehat \alpha \widehat {\mathfrak B} U^{110} + \widehat T \widehat \alpha U^{011} 
+ \widehat \alpha^2 U^{020} + \cdots.
$$
Replacing $\lambda _{0}$ by (\ref{lambda0}) in (\ref{exp lambda b}) and
identifying monomials leads to%
\begin{eqnarray*}
ia_{1}H_{0}U^{0} &=&(L_{0}-i\omega _{0}H_{0})U^{010}+L_{010}U^{0} { - i \omega_0 H_{010} U^0}, \\
a_{3}H_{0}U^{0} &=&(L_{0}-i\omega _{0}H_{0})U^{001}+L_{001}U^{0}, \\
ib_{1}H_{0}U^{0} &=&(L_{0}-i\omega _{0}H_{0})U^{100}+L_{100}U^{0}, \\
a_{4}H_{0}U^{0}+ib_{1}H_{0}U^{100} &=&(L_{0}-i\omega
_{0}H_{0})U^{200}+L_{100}U^{100}, \\
a_{2}H_{0}U^{0}+ia_{1}H_{0}U^{010} &=&(L_{0}-i\omega
_{0}H_{0})U^{020}+L_{010}U^{010}+L_{020}U^{0}
\\ & & 
{ - i a_1 H_{010} U^0 - i \omega_0 H_{010} U^{010}}
\end{eqnarray*}%
and
\begin{align*}
& 2a_{5}H_{0}U^{0}+ia_{1}H_{0}U^{100}+ib_{1}H_{0}U^{010}+i\omega_{0}H_{010}U^{100}  + i b_1 H_{010} U^0 
\\
&=(L_{0}-i\omega_{0}H_{0})U^{110}+L_{010}U^{100}+L_{100}U^{010}
+ L_{110} U^0,
\end{align*}%
where $U^{0}=(u_{x}^{0},u_{y}^{0})^{t}$ is the eigenvector satisfying%
\begin{equation*}
(i\omega _{0}H_{0}-L_{0})U^{0}=0.
\end{equation*}%
Let $V^0 = (v_x^0,v_y^0)$ be the adjoint vector, which satisfies
\begin{eqnarray}
(\alpha _{c}^{2}-D^{2})(-i\omega _{0}+\alpha _{c}^{2}-D^{2}+i\mathfrak{B}%
_{c}x)v_{x}^{0}-v_{y}^{0} &=&0,  \label{adj eigenvect} \\
(-i\omega _{0}+\alpha _{c}^{2}-D^{2}+i\mathfrak{B}_{c}x)v_{y}^{0}-\alpha
_{c}^{2}T_{c}g(x)v_{x}^{0} &=&0. \notag
\end{eqnarray}%
Then taking the scalar product of the previous equation with $V^0$, we obtain
\begin{eqnarray*}
ia_{1}\langle H_{0}U^{0},V^{0}\rangle  &=&\langle L_{010}U^{0},V^{0}\rangle 
{ - i \omega_0 \langle H_{010} U^0, V^0 \rangle},
\\
a_{3}\langle H_{0}U^{0},V^{0}\rangle  &=&\langle L_{001}U^{0},V^{0}\rangle ,
\\
ib_{1}\langle H_{0}U^{0},V^{0}\rangle  &=&\langle L_{100}U^{0},V^{0}\rangle ,
\\
a_{4}\langle H_{0}U^{0},V^{0}\rangle  &=&-ib_{1}\langle
H_{0}U^{100},V^{0}\rangle +\langle L_{100}U^{100},V^{0}\rangle , \\
a_{2}\langle H_{0}U^{0},V^{0}\rangle  &=&-ia_{1}\langle
H_{0}U^{010},V^{0}\rangle +\langle L_{010}U^{010}+L_{020}U^{0},V^{0}\rangle 
\\
& & 
{ - i a_1 \langle H_{010} U^0,V^0 \rangle
- i \omega_0 \langle H_{010} U^{010}, V^0 \rangle},
\\
{ 2} a_{5}\langle H_{0}U^{0},V^{0}\rangle  &=&-i\langle
a_{1}H_{0}U^{100}+b_{1}H_{0}U^{010}+\omega _{0}H_{010}U^{100},V^{0}\rangle
\\
&&
+\langle L_{010}U^{100}+L_{100}U^{010},V^{0}\rangle
{ - i b_1 \langle H_{010} U^0, V^0 \rangle + \langle L_{110} U^0, V^0 \rangle},
\end{eqnarray*}%
where, by construction,  $a_{1}$ and $b_{1}$ are real. We observe that
\begin{equation*}
\langle H_{0}U^{0},V^{0}\rangle =\langle \zeta _{1},\zeta _{1}^{\ast
}\rangle 
\end{equation*}%
which is given by (\ref{scal prof dzeta1}).

Once the eigenvectors $U^{0}$ and $V^{0}$ are computed, 
we can compute $a_1$, $a_3$ and $b_1$. We can then evaluate $U^{010}$, $U^{001}$ and $U^{100}$,
which leads to the computation of $a_2$, $a_4$ and $a_5$.

\medskip

Numerically, we check that, as expected, $a_{3r}>0,$ $a_{4r}<0,$ $a_{2r}<0$ and that
{ $a_{4r}a_{2r} > a_{5r}^{2}$.}

%%%%%%%%%%%%%%%%%%%%%%%%%%%%%%%%%%%%%%%%%%%%%%%%%%%%%%%%%%%%%%%%%%%%%%%%%

\subsection{Computation of $b$ and $c$ \label{appendixB}}

%%%%%%%%%%%%%%%%%%%%%%%%%%%%%%%%%%%%%%%%%%%%%%%%%%%%%%%%%%%%%%%%%%%%%%%%%

We start \ with the limit Navier-Stokes system (\ref{NSlimit}), where we look
for solutions which are $2\pi /\beta _{c}$ periodic in $y.$ Below, we suppress the hat on $u_{y}$. Let us define the operator $%
L_{0}$ as
\begin{equation}
L_{0}U=\Delta _{\bot }U+x\mathfrak{R}\partial _{y}U-\nabla _{\bot
}p_{0}+
{  \left(  \begin{array}{c}
T_{c}g(x)u_{y} \\ 
0 \\ 
u_{x}%
\end{array}%
\right) },  \label{defL0 a}
\end{equation}%
with%
\begin{equation}
\nabla _{\bot }\cdot U_{\bot }+\mathfrak{R}\partial _{y}u_{y}=0,
\label{defL0 b}
\end{equation}%
where $p_{0}$ is such that $L_{0}U$ satisfies the boundary conditions
\begin{equation}
L_{0}U|_{x} =0, \quad \hbox{when} \quad x=\pm 1/2,  \label{def L0 c} 
\end{equation}
and the modified divergence free condition
$$
\nabla _{\bot }\cdot (L_{0}U)_{\bot }+\mathfrak{R}\partial _{y}(L_{0}U)_{y} = 0.  
$$
Now we define the symmetric quadratic operator $B$ as%
\begin{equation}
B(U,V)=-\frac{1}{2}[(U_{\bot }\cdot \nabla _{\bot })V+(V_{\bot }\cdot \nabla
_{\bot })U]+\frac{T(1-\mu )}{2}\left( 
\begin{array}{c}
u_{y}v_{y} \\ 
0 \\ 
0%
\end{array}%
\right) +\nabla _{\bot }q,  \label{def B}
\end{equation}%
with the boundary conditions
$$
B(U,V)|_{x}  =0, \quad \hbox{when} \quad x=\pm 1/2, 
$$
and the modified divergence free condition
$$
\nabla _{\bot }\cdot (B(U,V))_{\bot }+\mathfrak{R}\partial _{y}B(U,V)_{y}=0.
$$
Now, our system for $T=T_{c}$ may be written as%
\begin{equation}
\partial _{t}U=L_{0}U+B_{c}(U,U),  \label{reduced NS}
\end{equation}%
where 
\begin{equation*}
B_{c}=B|_{T_{c}}.
\end{equation*}%
To obtain the coefficients $b$ and $c$ of the amplitude equations, we
first decompose $U$ as%
\begin{equation}
U=A\zeta _{1}+B\zeta _{2}+\overline{A}\overline{\zeta _{1}}+\overline{B}%
\overline{\zeta _{2}}+\Phi (A,\overline{A},B,\overline{B}),
\label{centermanifold}
\end{equation}%
where we expand $\Phi $ as%
\begin{equation*}
\Phi =\sum\limits_{p+q+r+s\geq 2}\Phi _{pqrs}A^{p}\overline{A}^{q}B^{r}%
\overline{B}^{s},
\end{equation*}%
\begin{eqnarray*}
\zeta _{1} &=&e^{i(\alpha _{c}z+\beta _{c}y)}U_{1}(x),\zeta
_{2}=e^{i(-\alpha _{c}z+\beta _{c}y)}SU_{1}(x), \\
(L_{0}-i\omega _{0})\zeta _{1} &=&0,(L_{0}-i\omega _{0})\zeta _{2}=0,
\end{eqnarray*}%
where coefficients $\Phi _{pqrs}$ are divergence free and satisfy the same
boundary conditions as $U.$ Now the principal part of the amplitude equations
reduces to 
\begin{eqnarray}
\frac{dA}{dt} &=&i\omega _{0}A+bA|A|^{2}+cA|B|^{2},
\label{reducedamplitudes} \\
\frac{dB}{dt} &=&i\omega _{0}B+cB|A|^{2}+bB|B|^{2},  \notag
\end{eqnarray}%
and the coefficients $b$ and $c$ are deduced from an identification of the
monomials, first in computing $\partial _{t}U$ in (\ref{centermanifold}) in
using (\ref{reducedamplitudes}), second in replacing $U$ in (\ref{reduced NS}%
). The identification of the monomials $A^{2},|A|^{2},|B|^{2},A|A|^{2},A|B|^{2}$
lead to%
\begin{eqnarray*}
0 &=&(L_{0}-2i\omega _{0})\Phi _{2000}+B(\zeta _{1},\zeta _{1}), \\
0 &=&L_{0}\Phi _{1100}+2B(\zeta _{1},\overline{\zeta _{1}}), \\
0 &=&L_{0}\Phi _{0011}+2B(\zeta _{2},\overline{\zeta _{2}}), \\
0 &=&(L_{0}-2i\omega _{0})\Phi _{1010}+2B(\zeta _{1},\zeta _{2}), \\
0 &=&L_{0}\Phi _{1001}+2B(\zeta _{1},\overline{\zeta _{2}}),
\end{eqnarray*}%
\begin{eqnarray*}
b\zeta _{1} &=&(L_{0}-i\omega _{0})\Phi _{2100}+2B(\zeta _{1},\Phi
_{1100})+2B(\overline{\zeta _{1}},\Phi _{2000}), \\
c\zeta _{1} &=&(L_{0}-i\omega _{0})\Phi _{1011}+2B(\zeta _{1},\Phi
_{0011})+2B(\overline{\zeta _{2}},\Phi _{1010})+2B(\zeta _{2},\Phi _{1001}),
\end{eqnarray*}%
where we notice that $L_{0}$ and $(L_{0}-2i\omega _{0})$ are invertible,
while $(L_{0}-i\omega _{0})$ has a two-dimensional kernel spanned by $(\zeta
_{1},\zeta _{2})$, so that Fredholm alternative gives coefficients $b$ and $c
$ in deriving the compatibility condition for each of the two last equations. 
Using the results of section \ref{appendixrange}, we obtain that
the coefficients $b$ and $c$ are determined by
\begin{eqnarray}
b\langle \zeta _{1},\zeta _{1}^{\ast }\rangle &=&\langle 2B(\zeta _{1},\Phi
_{1100})+2B(\overline{\zeta _{1}},\Phi _{2000}),\zeta _{1}^{\ast }\rangle ,
\label{coef b} \\
c\langle \zeta _{1},\zeta _{1}^{\ast }\rangle &=&\langle 2B(\zeta _{1},\Phi
_{0011})+2B(\overline{\zeta _{2}},\Phi _{1010})+2B(\zeta _{2},\Phi
_{1001}),\zeta _{1}^{\ast }\rangle ,  \label{coef c}
\end{eqnarray}%
with%
\begin{eqnarray*}
\Phi _{2000} &=&(2i\omega _{0}-L_{0})^{-1}B(\zeta _{1},\zeta _{1}), \\
\Phi _{1100} &=&-2L_{0}^{-1}B(\zeta _{1},\overline{\zeta _{1}}), \\
\Phi _{0011} &=&-2L_{0}^{-1}B(\zeta _{2},\overline{\zeta _{2}})=S\Phi
_{1100}, \\
\Phi _{1010} &=&2(2i\omega _{0}-L_{0})^{-1}B(\zeta _{1},\zeta _{2}), \\
\Phi _{1001} &=&-2L_{0}^{-1}B(\zeta _{1},\overline{\zeta _{2}}).
\end{eqnarray*}

%%%%%%%%%%%%%%%%%%%%%%%%%%%%%%%%%%%%%%%%%%%%%%%%%%%%

\subsubsection{Computation of the various $\Phi$}

%%%%%%%%%%%%%%%%%%%%%%%%%%%%%%%%%%%%%%%%%%%%%%%%%%%%

We now successively compute $\Phi_{2000}$, $\Phi_{1100}$, $\Phi_{0011}$, $\Phi_{1010}$, $\Phi_{1001}$.
We have%
\begin{equation*}
B(\zeta _{1},\zeta _{1})=e^{2i(\alpha _{c}z+\beta
_{c}y)}(B_{x}^{2000},B_{y}^{2000},B_{z}^{2000})^t +e^{2i(\alpha _{c}z+\beta
_{c}y)}(Dp,0,2i\alpha _{c}p)^t,
\end{equation*}%
with%
\begin{equation*}
B_{x}^{2000}=i\mathfrak{B}_{c}u_{x}^{0}u_{y}^{0}+T_{c}\frac{(1-\mu )}{2}%
(u_{y}^{0})^{2},
\end{equation*}%
\begin{equation*}
B_{y}^{2000}=-u_{x}^{0}Du_{y}^{0}+u_{y}^{0}Du_{x}^{0}+i\mathfrak{B}%
_{c}(u_{y}^{0})^{2},
\end{equation*}%
\begin{equation*}
B_{z}^{2000} = \frac{i}{\alpha_{c}} \Big[ (D u_{x}^{0})^{2} - u_{x}^{0} D^{2} u_{x}^{0} - \mathfrak{B}_{c}^{2} (u_{y}^{0})^{2} \Big]
+ \frac{\mathfrak{B}_{c}}{\alpha_{c}} \Big( u_{x}^{0} D u_{y}^{0} - 2 D u_{x}^{0} u_{y}^{0} \Big).
\end{equation*}%
Thus
\begin{equation*}
\Phi _{2000}=e^{2i(\alpha _{c}z+\beta _{c}y)}(\phi _{x}^{2000},\phi
_{y}^{2000},\phi _{z}^{2000})^{t},
\end{equation*}%
with%
\begin{eqnarray*}
&&(2i\omega _{0}+4\alpha _{c}^{2}-D^{2}-2i\mathfrak{B}_{c}x)(4\alpha
^{2}-D^{2})\phi _{x}^{2000}-4\alpha ^{2}T_{c}g(x)\phi _{y}^{2000} \\
&& \qquad = 4\alpha ^{2}B_{x}^{2000}+2i\mathfrak{B}_{c}DB_{y}^{2000}+2i\alpha
_{c}DB_{z}^{2000},
\end{eqnarray*}%
\begin{equation*}
(2i\omega _{0}+4\alpha _{c}^{2}-D^{2}-2i\mathfrak{B}_{c}x)\phi
_{y}^{2000}-\phi _{x}^{2000}=B_{y}^{2000},
\end{equation*}%
\begin{equation*}
\phi _{z}^{2000}=\frac{i}{2\alpha _{c}}(D\phi _{x}^{2000}+2i\mathfrak{B}%
_{c}\phi _{y}^{2000}),
\end{equation*}%
\begin{equation*}
\phi _{x}^{2000}=D\phi _{x}^{2000}=\phi _{y}^{2000}=0,x=\pm 1/2.
\end{equation*}
We now turn to the computation of $\Phi _{1100}$.
We have%
\begin{equation*}
2B(\zeta _{1},\overline{\zeta _{1}}%
)=(B_{x}^{1100},B_{y}^{1100},B_{z}^{110})^{t}+(Dp,0,0)^{t},
\end{equation*}%
with%
\begin{equation*}
B_{x}^{1100}=-2D(u_{x}^{0}\overline{u_{x}^{0}})+i\mathfrak{B}_{c}(u_{x}^{0}%
\overline{u_{y}^{0}}-\overline{u_{x}^{0}}u_{y}^{0})+T_{c}(1-\mu )u_{y}^{0}%
\overline{u_{y}^{0}},
\end{equation*}%
\begin{equation*}
B_{y}^{1100}=-D(u_{x}^{0}\overline{u_{y}^{0}}+\overline{u_{x}^{0}}u_{y}^{0}),
\end{equation*}%
\begin{equation*}
B_{z}^{1100}=\frac{i}{\alpha _{c}}(u_{x}^{0}D^{2}\overline{u_{x}^{0}} -%
\overline{u_{x}^{0}}D^{2}u_{x}^{0})+\frac{\mathfrak{B}_{c}}{\alpha _{c}}%
(u_{x}^{0}D\overline{u_{y}^{0}} +  \overline{u_{x}^{0}}Du_{y}^{0}).
\end{equation*}%
Thus
\begin{equation*}
\Phi _{1100}=(0,\phi _{y}^{1100},\phi _{z}^{1100})^{t},
\end{equation*}%
with%
\begin{eqnarray*}
-T_{c}g(x)\phi _{y}^{1100} &=&B_{x}^{1100} { -Dp}, \\
-D^{2}\phi _{y}^{1100} &=&B_{y}^{1100}, \\
-D^{2}\phi _{z}^{1100} &=&B_{z}^{1100},
\end{eqnarray*}%
\begin{equation*}
\phi _{y}^{1100}=\phi _{z}^{1100}=0,x=\pm 1/2.
\end{equation*}%
We observe that $\phi _{y}^{1100}$ is real while $\phi _{z}^{1100}$ is pure
imaginary.
We now compute $\Phi _{0011}$.
We have 
\begin{equation*}
\Phi _{0011}=S\Phi _{1100},
\end{equation*}%
hence%
\begin{equation*}
\Phi _{0011}=(0,\phi _{y}^{1100},-\phi _{z}^{1100})^{t}=\overline{\Phi
_{1100}}.
\end{equation*}
Next to compute $\Phi _{1010}$, we observe that
\begin{equation*}
2B(\zeta _{1},\zeta _{2})=e^{2i\beta
_{c}y}(B_{x}^{1010},B_{y}^{1010},B_{z}^{1010})^{t}+e^{2i\beta
_{c}y}(Dp,0,0)^{t},
\end{equation*}%
with%
\begin{equation*}
B_{x}^{1010}=-4u_{x}^{0}Du_{x}^{0}-2i\mathfrak{B}%
_{c}u_{x}^{0}u_{y}^{0}+T_{c}(1-\mu )(u_{y}^{0})^{2},
\end{equation*}%
\begin{equation*}
 B_{y}^{1010}= - 2(u_{y}^{0}Du_{x}^{0}+u_{x}^{0}Du_{y}^{0})
-2i\mathfrak{B}_{c}(u_{y}^{0})^{2} ,
\end{equation*}%
\begin{equation*}
B_{z}^{1010}=0.
\end{equation*}%
Thus
\begin{equation*}
\Phi _{1010}=e^{2i\beta _{c}y}(\phi _{x}^{1010},\phi _{y}^{1010},0)^{t},
\end{equation*}%
with%
%\begin{equation*}
%(2i\omega _{0}-D^{2}-2i\mathfrak{B}_{c}x)\phi _{x}^{1010}-T_{c}g(x)\phi
%_{y}^{1010}=B_{x}^{1010} {\color{red} +Dp},
%\end{equation*}%
\begin{equation*}
(2i\omega _{0}-D^{2}-2i\mathfrak{B}_{c}x)D\phi _{x}^{1010}+2i\mathfrak{B}%
_{c}\phi _{x}^{1010}=  -2i\mathfrak{B}_{c}B_{y}^{1010},
\end{equation*}%
\begin{equation*}
D\phi _{x}^{1010}+2i\mathfrak{B}_{c}\phi _{y}^{1010}=0,
\end{equation*}%
\begin{equation*}
\phi _{x}^{1010}=D\phi _{x}^{1010}=\phi _{y}^{1010}=0,x=\pm 1/2.
\end{equation*}
It remains to compute $\Phi _{1001}$.
Note that 
\begin{equation*}
2B(\zeta _{1},\overline{\zeta _{2}})=e^{2i\alpha
_{c}z}(B_{x}^{1001},B_{y}^{1001},B_{z}^{1001})^{t}+e^{2i\alpha
_{c}z}(Dp,0,2i\alpha _{c}p)^{t},
\end{equation*}%
with%
\begin{equation*}
B_{x}^{1001}=i\mathfrak{B}_{c}(\overline{u_{x}^{0}}u_{y}^{0}-u_{x}^{0}%
\overline{u_{y}^{0}})+T_{c}(1-\mu )u_{y}^{0}\overline{u_{y}^{0}},
\end{equation*}%
\begin{equation*}
B_{y}^{1001}=u_{y}^{0}D\overline{u_{x}^{0}}+\overline{u_{y}^{0}}%
Du_{x}^{0}-u_{x}^{0}D\overline{u_{y}^{0}}-\overline{u_{x}^{0}}Du_{y}^{0},
\end{equation*}%
\begin{eqnarray*}
B_{z}^{1001} &=&-\frac{i}{\alpha _{c}}[u_{x}^{0}D^{2}\overline{u_{x}^{0}}+%
\overline{u_{x}^{0}}D^{2}u_{x}^{0}-2Du_{x}^{0}D\overline{u_{x}^{0}}-2%
\mathfrak{B}_{c}^{2}u_{y}^{0}\overline{u_{y}^{0}}] \\
&&-\frac{\mathfrak{B}_{c}}{\alpha _{c}}[u_{x}^{0}D\overline{u_{y}^{0}}-%
\overline{u_{x}^{0}}Du_{y}^{0}+2u_{y}^{0}D\overline{u_{x}^{0}}-2\overline{%
u_{y}^{0}}Du_{x}^{0}].
\end{eqnarray*}%
We notice that $B_{x}^{1001}$ and $B_{y}^{1001}$ are real, while $%
B_{z}^{1001}$ is imaginary. Then we get
\begin{equation*}
\Phi ^{1001}=e^{2i\alpha _{c}z}(\phi _{x}^{1001},\phi _{y}^{1001},\phi
_{z}^{1001})^{t},
\end{equation*}%
with%
\begin{equation*}
(4\alpha _{c}^{2}-D^{2})^{2}\phi _{x}^{1001}-4\alpha _{c}^{2}T_{c}g(x)\phi
_{y}^{1001}=4\alpha _{c}^{2}B_{x}^{1001}+2i\alpha _{c}DB_{z}^{1001},
\end{equation*}%
\begin{equation*}
(4\alpha _{c}^{2}-D^{2})\phi _{y}^{1001}-\phi _{x}^{1001}=B_{y}^{1001},
\end{equation*}%
\begin{equation*}
\phi _{z}^{1001}=\frac{i}{2\alpha _{c}}D\phi _{x}^{1001},
\end{equation*}%
\begin{equation*}
\phi _{x}^{1001}=D\phi _{x}^{1001}=\phi _{y}^{1001}=0,x=\pm 1/2.
\end{equation*}%
We notice that $\phi _{x}^{1001}$ and $\phi _{y}^{1001}$ are real, while $%
\phi _{z}^{1001}$ is imaginary (this is due to symmetry $S).$

%%%%%%%%%%%%%%%%%%%%%%%%%%%%%%%%%%%%%%%%%%%%%%%%%%%%%%%%%%%%%%%%%

\subsubsection{Computation of the terms involving $B$}

%%%%%%%%%%%%%%%%%%%%%%%%%%%%%%%%%%%%%%%%%%%%%%%%%%%%%%%%%%%%%%%%%

We have%
\begin{equation*}
2B(\zeta _{1},\Phi _{1100})=e^{i(\alpha _{c}z+\beta
_{c}y)}(B_{x}^{1,1100},B_{y}^{1,1100},B_{z}^{1,1100})^{t}+e^{i(\alpha
_{c}z+\beta _{c}y)}(Dp_{1},0,i\alpha _{c}p_{1})^{t},
\end{equation*}%
with%
\begin{eqnarray*}
B_{x}^{1,1100} &=&-i\alpha _{c}u_{x}^{0}\phi _{z}^{1100}+T_{c}(1-\mu
)u_{y}^{0}\phi _{y}^{1100}, \\
B_{y}^{1,1100} &=&-u_{x}^{0}D\phi _{y}^{1100}-i\alpha _{c}u_{y}^{0}\phi
_{z}^{1100}, \\
 B_{z}^{1,1100} &=&-D(u_{x}^{0}\phi _{z}^{1100})-i\mathfrak{B}%
_{c}u_{y}^{0}\phi _{z}^{1100}.
\end{eqnarray*}
Next,
\begin{equation*}
2B(\overline{\zeta _{1}},\Phi _{2000})=e^{i(\alpha _{c}z+\beta
_{c}y)}(B_{x}^{\overline{1},2000},B_{y}^{\overline{1},2000},B_{z}^{\overline{%
1},2000})^{t}+e^{i(\alpha _{c}z+\beta _{c}y)}(Dp_{1},0,i\alpha
_{c}p_{1})^{t},
\end{equation*}%
with%
\begin{eqnarray*}
B_{x}^{\overline{1},2000} &=&-\overline{u_{x}^{0}}D\phi _{x}^{2000}-3D%
\overline{u_{x}^{0}}\phi _{x}^{2000}+2i\mathfrak{B}_{c}\overline{u_{y}^{0}}%
\phi _{x}^{2000}+i\alpha _{c}\overline{u_{x}^{0}}\phi
_{z}^{2000}+T_{c}(1-\mu )\overline{u_{y}^{0}}\phi _{y}^{2000}, \\
B_{y}^{\overline{1},2000} &=&-\overline{u_{x}^{0}}D\phi _{y}^{2000}-2D%
\overline{u_{x}^{0}}\phi _{y}^{2000}+2i\mathfrak{B}_{c}\overline{u_{y}^{0}}%
\phi _{y}^{2000}-D\overline{u_{y}^{0}}\phi _{x}^{2000}+i\alpha _{c}\overline{%
u_{y}^{0}}\phi _{z}^{2000}, \\
B_{z}^{\overline{1},2000} &=&D\overline{u_{x}^{0}}\phi _{z}^{2000}-\overline{%
u_{x}^{0}}D\phi _{z}^{2000}-i\mathfrak{B}_{c}\overline{u_{y}^{0}}\phi
_{z}^{2000}+\frac{i}{\alpha _{c}}D^{2}\overline{u_{x}^{0}}\phi _{x}^{2000}+%
\frac{\mathfrak{B}}{\alpha _{c}}D\overline{u_{y}^{0}}\phi _{x}^{2000}.
\end{eqnarray*}
We also have
\begin{equation*}
2B(\zeta _{1},\Phi _{0011})=e^{i(\alpha _{c}z+\beta
_{c}y)}(B_{x}^{1,0011},B_{y}^{1,0011},B_{z}^{1,0011})^{t}+e^{i(\alpha
_{c}z+\beta _{c}y)}(Dp_{2},0,i\alpha _{c}p_{2})^{t},
\end{equation*}%
with%
\begin{eqnarray*}
B_{x}^{1,0011} &=&i\alpha _{c}u_{x}^{0}\phi _{z}^{1100}+T_{c}(1-\mu
)u_{y}^{0}\phi _{y}^{1100}, \\
B_{y}^{1,0011} &=&-u_{x}^{0}D\phi _{y}^{1100}+i\alpha _{c}u_{y}^{0}\phi
_{z}^{1100}, \\
B_{z}^{1,0011} &=&D(u_{x}^{0}\phi _{z}^{1100})+i\mathfrak{B}%
_{c}u_{y}^{0}\phi _{z}^{1100}.
\end{eqnarray*}
Moreover,
\begin{equation*}
2B(\overline{\zeta _{2}},\Phi _{1010})=e^{i(\alpha _{c}z+\beta
_{c}y)}(B_{x}^{\overline{2},1010},B_{y}^{\overline{2},1010},B_{z}^{\overline{%
2},1010})^{t}+e^{i(\alpha _{c}z+\beta _{c}y)}(Dp_{3},0,i\alpha
_{c}p_{3})^{t},
\end{equation*}%
with%
\begin{eqnarray*}
B_{x}^{\overline{2},1010} &=&-D(\overline{u_{x}^{0}}\phi
_{x}^{1010})+T_{c}(1-\mu )\overline{u_{y}^{0}}\phi _{y}^{1010}, \\
B_{y}^{\overline{2},1010} &=&-\overline{u_{x}^{0}}D\phi _{y}^{1010}-D%
\overline{u_{y}^{0}}\phi _{x}^{1010}, \\
B_{z}^{\overline{2},1010} &=&(-\frac{i}{\alpha _{c}}D^{2}\overline{u_{x}^{0}}%
+\frac{\mathfrak{B}_{c}}{\alpha _{c}}\overline{Du_{y}^{0}})\phi _{x}^{1010}.
\end{eqnarray*}
The last computation is
\begin{equation*}
2B(\zeta _{2},\Phi _{1001})=e^{i(\alpha _{c}z+\beta
_{c}y)}(B_{x}^{2,1001},B_{y}^{2,1001},B_{z}^{2,1001})^{t}+e^{i(\alpha
_{c}z+\beta _{c}y)}(Dp_{4},0,i\alpha _{c}p_{4})^{t},
\end{equation*}%
with%
\begin{eqnarray*}
B_{x}^{2,1001} &=&-u_{x}^{0}D\phi _{x}^{1001}-3Du_{x}^{0}\phi _{x}^{1001}-2i%
\mathfrak{B}_{c}u_{y}^{0}\phi _{x}^{1001}+i\alpha _{c}u_{x}^{0}\phi
_{z}^{1001}+T_{c}(1-\mu )u_{y}^{0}\phi _{y}^{1001}, \\
B_{y}^{2,1001} &=&-u_{x}^{0}D\phi _{y}^{1001}-2Du_{x}^{0}\phi
_{y}^{1001}-Du_{y}^{0}\phi _{x}^{1001}-2i\mathfrak{B}_{c}u_{y}^{0}\phi
_{y}^{1001}+i\alpha _{c}u_{y}^{0}\phi _{z}^{1001}, \\
B_{z}^{2,1001} &=&-D(u_{x}^{0}\phi _{z}^{1001})+\frac{i}{\alpha _{c}}%
D^{2}u_{x}^{0}\phi _{x}^{1001}-\frac{\mathfrak{B}_{c}}{\alpha _{c}}%
Du_{y}^{0}\phi _{x}^{1001}-i\mathfrak{B}_{c}u_{y}^{0}\phi _{z}^{1001}.
\end{eqnarray*}

%%%%%%%%%%%%%%%%%%%%%%%%%%%%%%%%%%%%%%%%%%%%%%%%%%%%%%%%%%%%%

\subsubsection{Study of $b$ and $c$}

%%%%%%%%%%%%%%%%%%%%%%%%%%%%%%%%%%%%%%%%%%%%%%%%%%%%%%%%%%%%%

We finally give the explicit expressions of $b$ and $c$.
Using (\ref{coef b}), (\ref{coef c}) and (\ref{scal prof dzeta1}), we obtain%
\begin{align}
b\langle \zeta _{1},\zeta _{1}^{\ast }\rangle &
= \int_{-1/2}^{1/2} \Bigl[ (B_{x}^{1,1100}+B_{x}^{\overline{1},2000})(\alpha
_{c}^{2}-D^{2})\overline{v_{x}^{0}}+(B_{y}^{1,1100}+B_{y}^{\overline{1}%
,2000})\overline{v_{y}^{0}} \Bigr] \, dx,  \label{formula for b,c} \\
c\langle \zeta _{1},\zeta _{1}^{\ast }\rangle &
=\int_{-1/2}^{1/2} \Bigl[ (B_{x}^{1,0011}+B_{x}^{\overline{2},1010}+B_{x}^{2,1001})(\alpha _{c}^{2}-D^{2})\overline{v_{x}^{0}}%
\\
&\qquad \qquad +(B_{y}^{1,0011}+B_{y}^{\overline{2},1010}+B_{y}^{2,1001})\overline{v_{y}^{0}} \Bigr] \, dx,  \notag
\end{align}%
with
\begin{equation*}
\langle \zeta _{1},\zeta _{1}^{\ast }\rangle 
=\int_{-1/2}^{1/2} \Bigl[ Du_{x}^{0}D\overline{v_{x}^{0}}+\alpha _{c}^{2}u_{x}^{0}\overline{v_{x}^{0}}+u_{y}^{0} \overline{v_{y}^{0}} \Bigr] \,dx,
\end{equation*}%
where $u_{x}^{0},u_{y}^{0},v_{x}^{0},v_{y}^{0}$ are defined by (\ref{eigenvect}) and (\ref{adj eigenvect}).

\medskip 

The numerical evaluation of the coefficients $b$ and $c$ follows the algebraic procedure detailed in this appendix. The interval $x\in[-1/2,1/2]$ is discretised with $N$ Chebyshev collocation points and the differential operators $D_x$ and $D_x^2$ are replaced by their spectral differentiation matrices. The eigenvalue problem \eqref{eigenvect} is solved together with the boundary conditions to obtain the direct eigenvector $(u_x^0,u_y^0)$. The adjoint eigenvector $(v_x^0,v_y^0)$ is obtained from the discretised adjoint problem \eqref{adj eigenvect}, and its normalisation is chosen so that the scalar product $\langle\zeta_1,\zeta_1^*\rangle$ equals unity. The linear systems that define the second‑order harmonics $\Phi_{2000},\Phi_{1100},\Phi_{0011},\Phi_{1010},\Phi_{1001}$ are then assembled. Except for the Poisson problems for $\Phi_{1100}$ and $\Phi_{1001}$, each system is of the same structure as \eqref{eigenvect} and is solved with the appropriate factor $(\alpha,\mathfrak{B}_c,\omega_0)$. The quadratic interaction terms are built from the discretised eigenfunctions and harmonics, projected onto the adjoint eigenvector, and integrated using Clenshaw-Curtis quadrature. The values of $b$ and $c$ are finally obtained from formulas \eqref{formula for b,c}. All computations are performed with extended precision to guarantee well‑converged results.

\medskip

For $\mu = -1$, we obtain
$$
b =   -0.3931 + 1.6743i, \qquad c =       -0.3846 - 0.5753i.
$$

\subsubsection*{Acknowledgment}
D. Bian is partly supported by NSF of China under the contract 12271032.

%%%%%%%%%%%%%%%%%%%%%%%%%%%%%%%%%%%%%%%%%%%%%%%%%%%%%%%%%%%%%%%%%%%%%%%%%%%%%%%%%%%%%

\end{document}